\input amstex
\documentstyle{amsppt} 




\def\cal{\Cal}
\define\Y{\Cal Y}
\define\X{\Cal X}
\define\tlY{\Y_{-2}}
\define\tlX{\X}
\define\M{\Cal M}
\define\V{\Cal V}
\define\F{\Cal F}
\define\W{\Cal W}

\define\R{{\Bbb R}}
\define\Q{{\Bbb Q}}
\define\Z{{\Bbb Z}}
\define\C{{\Bbb C}}

\define\a{\alpha}
\undefine\b
\define\b{\beta}
\define\g{\gamma}
\undefine\l
\define\l{\lambda}

\define\r{\rho}

\define\lam{\Lambda^{-1,-1}}
\define\nn{\eta_{-}}
\define\np{\eta_{+}}
\define\nz{\eta_0}
\define\pip{\pi_{+}}
\define\pin{\pi_{-}}
\define\piz{\pi_0}
\define\pil{\pi_{\Lambda}}
\define\piq{\pi_{\frak t}}

\define\G{\text{\rm G}}
\define\Gc{\text{\rm G}_{\Bbb C}}
\define\tlG{\text{\rm H}}
\undefine\gg
\define\gg{\frak g}
\define\hh{\frak h}
\undefine\aa
\define\aa{\frak a}
\define\zz{\frak z}
\define\xx{x}
\define\n{\eurm n}
\define\h{\eurm h}
\define\UU{\frak U}

\define\Span{\operatorname{span}}
\define\Hom{\operatorname{Hom}}

\define\half{\frac{1}{2}}
\define\hph{\hphantom}
\define\pd{\partial}
\define\rel{{}^r}
\undefine\Im
\define\Im{\text{\rom{Im}}}
\define\ad{\text{\rom{ad}}\,}
\define\Ad{\text{\rom{Ad}}}
\define\LL{\text{\rom{L}}}


\TagsOnRight
\TagsAsMath
\loadmsbm
\loadeurm
\loadeusm

\topmatter

\abstract We prove an analog of Schmid's $\text{\rm SL}_2$-orbit theorem
for a class of variations of mixed Hodge structure which includes logarithmic
deformations, degenerations of 1-motives and archimedean heights.  In 
particular, as consequence this theorem, we obtain a simple formula for
the asymptotic behavior of the archimedean height of a flat family of 
algebraic cycles which depends only on the weight filtration and local
monodromy.
\endabstract

\author Gregory Pearlstein
\endauthor


\title $\text{\rm SL}_2$-Orbits and Degenerations of Mixed Hodge
        Structure 
\endtitle

\thanks I wish to thank the Max--Planck Institut f\"ur Mathematik, Bonn 
and the University of Massachusetts, Amherst for their generous hospitality 
during the preparation of this manuscript.  I also wish to thank Richard Hain 
for suggesting the applications of this work to Arakelov geometry presented 
in \S 5.
\endthanks

\endtopmatter

\document


\head \S 1.\quad Introduction \endhead

\par Let $f:X\to S$ be a smooth, projective morphism of complex,
quasi-projective varieties.  Then, by the work of Griffiths [G],
the cohomology groups $\V_s = H^k(X_s)$ patch together to 
form a variation of Hodge structure $\V$ over $S$.  Furthermore, 
as a consequence of Schmid's orbit theorems [S], [CKS], one has a 
complete local theory regarding how such variations of Hodge structure 
degenerate along the boundary of a (partial) compactification 
$S\hookrightarrow\bar S$. 

\par Namely, by the work of Hironaka [Hiro] and Borel [D1], we can restrict
our attention to the case where $S$ is a product of punctured disks
$\Delta^{*n}$ and the monodromy representation of $\V$ is given a by 
a system of unipotent transformations $T_j = e^{-N_j}$.  Schmid's 
nilpotent orbit theorem asserts that, after lifting the period
map of $\V$ to a $\pi_1$-equivariant map 
$$
	F:U^n\to\cal D						
$$
from a product of upper half-planes into the corresponding classifying
space of polarized Hodge structure, there exists an associated
nilpotent orbit
$$
     \theta(\text{\bf z}) = \exp(\sum_j\, z_j N_j).F_{\infty}
$$
which is asymptotic to $F(\text{\bf z})$ with respect to a suitable metric
on $\Cal D$.  

\par The possible nilpotent orbits $\theta(\text{\bf z})$ which can
arise in this way are, in turn, classified by the $\text{\rm SL}_2$-orbit 
theorem [S] [CKS] which, roughly speaking, says that every such nilpotent 
orbit $\theta(\text{\bf z})$ is asymptotic to another nilpotent orbit 
$\hat\theta(\text{\bf z})$ which arises from a representation
of $\text{\rm SL}_2(\R)^n$.  

\par More precisely, since the classifying space $\cal D$ is the homogeneous 
space of a real, semisimple Lie group $\G_{\R}$, one defines a 1-variable 
$\text{\rm SL}_2(\R)$ orbit to be a nilpotent orbit $\theta(z)$ for which 
there exists a base point $F_o\in\cal D$ and a Lie homomorphism 
$\psi:\text{\rm SL}_2(\R)\to\G_{\R}$ such that 
$$
       \theta(g.\sqrt{-1}) = \psi(g).F_o                           
$$
Schmid's 1-variable $\text{\rm SL}_2$-orbit theorem then asserts that given 
any nilpotent orbit $e^{zN}.F$ of pure, polarized Hodge structure, there 
exists a $\text{\rm SL}_2$-orbit $e^{zN}.\hat F$, and a distinguished real 
analytic function 
$$
        g:(a,\infty)\to\G_{\R}
$$
such that
\roster
\item"(a)" $e^{iyN}.F = g(y)e^{iyN}.\hat F$;
\item"(b)" $g(y)$ and $g^{-1}(y)$ have convergent series expansions about
$\infty$ of the form $(1+\sum_{k=1}^{\infty}\, A_k y^k)$ with 
$A_k\in\ker(\ad\, N)^{k+1}$.
\endroster

\par In this article, we consider analogous questions for morphisms
$f:X\to S$ which are no longer necessarily proper or smooth.  In this
context, the variations of pure Hodge structure considered above are
replaced (cf\. \S 3) by variations of graded-polarized mixed Hodge
structure which are {\it admissible} in the sense of Steenbrink and 
Zucker [SZ].

\par In [P3], we proved that for admissible variations over a 1-dimensional 
base $S$, one has a corresponding nilpotent orbit theorem.  To state our
main result, we recall (cf\. \S 2) that the period map of a variation of
graded-polarized mixed Hodge structure takes values in the quotient of a
classifying space $\M$ of graded-polarized mixed Hodge structure upon
which a Lie group $\G$ acts transitively by automorphisms.  Furthermore
[KP], in this setting the natural analogs of the $\text{\rm SL}_2$-orbits 
considered above are admissible nilpotent orbits $e^{zN}.\hat F$ for which 
the associated limiting mixed Hodge structure (cf\. \S 3) is split over $\R$.

\par Accordingly, by virtue of the above remarks, it is natural to conjecture 
that given an admissible nilpotent orbit $e^{zN}.F$, there should exists a 
split orbit $e^{zN}.\hat F$ and a distinguished real analytic function 
$$
       g:(a,\infty)\to\G
$$
such that
\roster 
\item"(a)" $e^{iyN}.F = g(y)e^{iyN}.\hat F$;
\item"(b)" $g(\infty) := \lim_{y\to\infty}\, g(y)\in\ker(\ad\, N)$;
\item"(c)" $g^{-1}(\infty)g(y)$ and $g^{-1}(y)g(\infty)$ have convergent 
series expansions about $\infty$ of the form $(1+\sum_{k>0}\, A_k y^{-k})$ 
with $A_k\in\ker(\ad N)^{k+1}$.
\endroster

\par In \S 6--9, we prove the existence [Theorem $(4.2)$] of such a 
function $g(y)$ provided the Hodge numbers of the associated classifying 
space $\M$ belong to one of the following two subcases, each of which 
arises in a number of geometric settings (e.g\. 1-motives [D2], logarithmic
deformations [U]): 
\roster
\item"(I)"  $h^{p,q} = 0$ unless $p+q = k$, $k-1$;
\item"(II)" $h^{p,q} = 0$ unless $p+q = 2k-1$, or $(p,q) = (k,k)$, 
$(k-1,k-1)$.
\endroster

\par In particular (cf\. \S 5), as a consequence of the 
$\text{\rm SL}_2$-orbit theorem described above, we obtain a simple 
formula for the  asymptotic behavior of the archimedean height  [Arak] [Beil] 
[GS] 
$$
       h(s) = \langle Z_s, W_s\rangle_{\infty}                     
$$
of a flat family of algebraic cycles $Z_s$, $W_s\subseteq X_s$ over a smooth
curve $S$, which depends only on the weight filtration and local monodromy
of the associated variation of mixed Hodge structure [H].  

\par As in [S] [CKS], the proof of Theorem $(4.2)$ boils down to the 
construction of an explicit solution to an associated system of 
\lq\lq monopole equations\rq\rq{} attached to the nilpotent orbit $e^{zN}.F$. 
More precisely (cf\. \S 2), in each of the two subcases $(\text{I})$ and
$(\text{II})$ considered above, there exists a natural subgroup $\tlG$ of
$\G$ which acts transitively on the corresponding classifying space $\M$
by isometries.  As such (cf\. \S 6), each choice of base point $F_o\in\M$ 
defines an auxiliary principal bundle
$$
	\tlG^{F_o} \to \tlG \to \tlG/\tlG^{F_o}
$$
$P$ over $\M$.  Accordingly, a choice of connection $\nabla$ on $P$
determines a lift of $e^{iyN}$ to an $\tlG$-valued function $h(y)$ which is 
tangent to $\nabla$.  Moreover, as in [S], the resulting
function $h(y)$ satisfies a differential equation [Theorem $(6.11)$] of 
the form
$$
     h^{-1}\frac{dh}{dy} = -\LL\,\Ad(h^{-1}(y))N      \tag{1.1}
$$
relative to a suitable endomorphism $\LL$ of $\hh = Lie(\tlG)$.
In particular, as a consequence of equation $(1.1)$, the Hodge components
$$
        \b(y) = \b^{1,-1}(y) + \b^{0,0}(y) + \b^{-1,1}(y)
                + \b^{0,1}(y) + \b^{-1,0}(y)
$$
of the function $\b(y) = \Ad(h^{-1}(y))N$ associated to a nilpotent orbit 
$e^{iyN}.F$ of type $(\text{I})$ satisfy the following system of differential 
equations
$$
         \frac{d}{dy}\b_0(y) = -[\b_0(y),\LL \b_0(y)],\qquad
         \b_0(y) = \sum_{r+s=0}\, \b^{r,s}(y)                      \tag{1.2}
$$
$$
           \frac{d}{dy}\pmatrix \b^{-1,0} \\ \b^{0,-1} \endpmatrix
           = \sqrt{-1}
             \pmatrix   \ad\, \b^{0,0}  & -2\,\ad\, \b^{-1,1} \\
                     2\,\ad\, \b^{1,-1} & -\ad\, \b^{0,0}  \endpmatrix
             \pmatrix \b^{-1,0} \\ \b^{0,-1} \endpmatrix
                                                                    \tag{1.3}
$$

\par Following [S], we then observe that equation $(1.2)$ becomes equivalent 
to Nahm's equations [Hitch]
$$
\gathered
  -2\frac{d}{dy} X^+(y) = [Z(y),X^+(y)],\qquad
   2\frac{d}{dy} X^-(y) = [Z(y),X^-(y)]                             \\
   -\frac{d}{dy} Z(y)   = [X^+(y),X^-(y)]
\endgathered                                                        \tag{1.4}
$$
upon setting $X^+(y) = 2i\b^{1,-1}(y)$, $Z(y) = 2i\b^{0,0}(y)$ and 
$X^-(y) = -2i\b^{-1,1}(y)$.  Moreover, using the methods of [CKS], one can 
construct a series solution (cf\. \S 7) to equation $(1.4)$ in the form of 
a function
$$
    \Phi(y):(a,\infty)\to \Hom(sl_2(\C),\gg_{\C}),\qquad
    \Phi(y) = \sum_{n\geq 0} \Phi_n y^{-1-n/2}                      \tag{1.5}
$$
such that $X^-(y) = \Phi(y)\xx^-$, $Z(y) = \Phi(y)\zz$, and 
$X^+(y) = \Phi(y)\xx^+$ where:
$$
	\xx^- = \half\pmatrix 1 &  -i  \\
			      -i & -1 \endpmatrix,\qquad	
	\zz = \pmatrix  0 & -i \\
	                i & 0 \endpmatrix,\qquad
	\xx^+ = \half\pmatrix 1 &  i \\
			      i & -1 \endpmatrix		\tag{1.6}
$$

\par Building upon the series solution $(1.5)$,
we then construct a similar  series solution to $(1.3)$ in \S 8.  
Taken with equation $(1.1)$, such an series solution for $\b(y)$ then allows 
us to compute $h(y)$ modulo left multiplication by an element $h_o\in\tlG$.  
Imposing the boundary condition
$$
       \lim_{y\to\infty} e^{-iyN}h(y).F_o = F
$$	
then determines $h_o$.  Having computed $h(y)$, the
desired function $g(y)$ is then given by the formula
$$
         h(y) = g(y) y^{-H/2}					
$$
where $H = \Phi_0(\xx^+ + \xx^-)$.
\vskip 3pt

\par To illustrate how the $\text{\rm SL}_2$-orbit theorem described above 
works in the context of a geometric example, let $X$ be a compact Riemann 
surface and 
$$
	c_1 = c_{12} - c_{11},\qquad c_2 =  c_{22} - c_{21}	\tag{1.7}
$$
be a pair of disjoint 0-cycles on $X$.  Then (up to an additive constant),
there exists a unique harmonic function $f:X-|c_2|\to\R$ such that 
$$
	\Omega = \frac{1}{2\pi}(*df - i\,df)		\tag{1.8}
$$
is a holomorphic 1-form on $X-|c_2|$ with simple poles along 
$|c_2| = \{c_{22},c_{21}\}$ and residues
$$
	\text{Res}_{c_{22}}(\Omega) = \frac{1}{2\pi i},\qquad 
	\text{Res}_{c_{21}}(\Omega) = -\frac{1}{2\pi i}
$$
The archimedean height of $c_1$ and $c_2$ is then defined to be
$$
	\langle c_1,c_2 \rangle 
	= -2\pi\,\Im\left(\int_{c_{11}}^{c_{12}}\,\Omega\right)	\tag{1.9}
$$

\par To bring in the mixed Hodge structures, we now recall [D2] that the
elements of $H^1(X-|c_2|)$ can be decomposed according to (mixed) Hodge
type.  Furthermore, with respect to this decomposition, $\Omega$ generates 
the classes of type $(1,1)$.  As such, the integral $(1.9)$
can be viewed as a period of $H^1(X-|c_2|)$ with respect to $c_1$.  
Therefore, upon varying the triple $(X,c_1,c_2)$, the integral $(1.9)$ 
defines a \lq\lq period map\rq\rq{} whose asymptotic behavior is governed by 
Theorem $(4.2)$.  In particular [Theorem $(5.19)$], near a degenerate point 
$s=0$,
$$
	\langle c_1(s), c_2(s) \rangle \approx -\mu\log|s|
$$
where $\mu$ is a constant which depends only on the local monodromy of 
the associated variation of mixed Hodge structure.

\par More concretely, let $E\to\Delta^*$ be the family of elliptic curves 
$$
	E_s = \C/(\Z\oplus\tau(s)\Z)
$$
defined by the function $\tau(s) = \frac{1}{\pi i}\log(s)$ and 
$$
	h(s) = \langle e_3 - e_0, e_2 - e_1\rangle
$$ 
be the height function determined by the 2-torsion points
$$
	e_0 = 0,\quad e_1 = \half,\quad e_2 = \frac{\tau}{2},\quad
	e_3 = \half(1+\tau)			
$$
Then, a short calculation shows that
$$
	h(s) = -\log\left|\frac{\vartheta^2(e_2)}{\vartheta^2(e_1)}\right|
	       + \half\log|\exp(-2\pi i e_3)|			
$$
where $\vartheta$ is Riemann's theta function, and hence
$h(s)\approx -\frac{1}{2}\log|s|$ as $s\to 0$.
\vskip 3pt

\par To illustrate another application of the $\text{\rm SL}_2$-orbit theorem,
let 
$$
	F:U\to\M					
$$
be the period map of a non-constant, admissible variation of type 
$(\text{I})$.  Then, as a consequence of Theorem $(4.2)$, the holomorphic 
sectional curvature of $F(z)$ is negative, and bounded away from zero as 
$\text{Im}(z)\to\infty$ [Theorem $(4.9)$].

\par Heuristically, the proof of this fact boils down to replacing $F(z)$ by 
the corresponding split orbit $\hat\theta(z) = e^{zN}.\hat F$ and then noting
that split orbits of type $(\text{I})$ are actually $\text{\rm SL}_2$-orbits. 
More precisely, by virtue of the above remarks,
$$
	||F_*(d/dz)||_{F(z)} \approx 
	||\hat\theta_*(d/dz)||_{\hat\theta(z)}		
$$
Accordingly, since $\hat\theta(z)$ is a nilpotent orbit, 
$\hat\theta_*(\frac{d}{dz})$ is basically just $N$, and hence
(up to a constant scalar factor)
$$
	||F_*(d/dz)||_{F(z)} 
	\approx ||N||_{\hat\theta(z)}			
$$
Therefore (cf\. \S 2), since the real elements of $\G$ act on $\M$ by
isometries, it then follows that
$$
	||N||_{\hat\theta(z)} = ||N||_{e^{xN}e^{iyN}.\hat F}
	= ||N||_{e^{iyN}.\hat F}		
$$
Consequently, since $\hat\theta(z)$ is actually an $\text{\rm SL}_2$-orbit, 
$$
	e^{iyN}.\hat F = \exp(-\half\log(y)H)e^{iN}.\hat F
$$
where $H$ is real and $[H,N] = -2N$.  Thus,
$$
\aligned
	||F_*(d/dz)||_{F(z)} &\approx ||N||_{e^{iyN}.\hat F_{\infty}}  
	 = ||N||_{\exp(-\half\log(y)H)e^{iN}.\hat F_{\infty}}          \\
	&= ||\Ad(\exp(\half\log(y)H))N||_{e^{iN}.\hat F_{\infty}}     \\  
	&= (1/y)||N||_{e^{iN}.\hat F_{\infty}}
\endaligned
$$
and hence the pullback of the metric of $\M$ along $F$ is asymptotic to
a constant multiple of the Poincar\'e metric.

\head \S 2.\quad Preliminary Remarks \endhead

\par In this section, we recall the construction of the period map of a
variation of graded-polarized mixed Hodge structure, and discuss the
geometry of the associated classifying spaces of graded-polarized mixed
Hodge structure [K] [P2] [U].

\definition{Definition 2.1} Let $S$ be a complex manifold.  Then, a variation
of graded-polarized mixed Hodge structure over $S$ consists of the following
data:
\roster
\item A finite rank, $\Q$-local system $\V_{\Q}$ over $S$;
\item A rational, increasing filtration 
$\cdots\subseteq \W_k\subseteq\W_{k+1}\subseteq\cdots$ of 
$\V_{\C} = \V_{\Q}\otimes\C$ by sublocal systems;
\item A decreasing filtration
$\cdots\subseteq \F^p\subseteq\F^{p-1}\subseteq\cdots$
of $\V = \V_{\C}\otimes{\Cal O}_S$ by holomorphic subbundles;

\item A collection of non-degenerate bilinear forms 
$$
	Q_k:Gr^{\W}_k(\V_{\Q})\otimes Gr^{\W}_k(\V_{\Q})\to\Q
$$
of alternating parity $(-1)^k$;
\endroster 
subject to the following two conditions:
\roster
\item"(a)" $\F$ is horizontal with respect to the Gauss--Manin connection 
$\nabla$ of $\V$,
i.e. $\nabla(\F^p)\subseteq \F^{p-1}\otimes\Omega^1_S$;
\item"(b)" For each index $k$, $(Gr^{\W}_k(\V_{\Q}),\F Gr^{\W}_k,Q_k)$ is a 
variation of pure, polarized Hodge structure of weight $k$.
\endroster
\enddefinition

\par In analogy with the pure case [S], the isomorphism class of a variation 
of graded-polarized mixed Hodge structure $\V\to S$ is determined by its 
period map
$$
	\varphi:S\to\Gamma\backslash \M,\qquad \Gamma = \text{Image}(\rho)
								\tag{2.2}
$$
and its monodromy representation $\rho:\pi_1(S,s_0)\to GL(\V_{s_o})$ on a 
fixed reference fiber $V = \V_{s_o}$.  More precisely, let $W$ and 
$Q=\{Q_k\}$ denote the specialization of the weight filtration and 
graded-polarizations of $\V$ to $V$.  Define $X$ to be the flag variety 
consisting of all decreasing filtrations $F$ of $V$ such that 
$$
	\dim(F^p) = \text{rank}(\F^p)
$$
and let $\M$ denote the classifying space [P2] consisting of all 
filtrations $F\in X$ such that $(F,W)$ is a mixed Hodge structure which is 
graded-polarized by $Q$.  Then, the period map $(2.2)$ is obtained by simply
pulling back the Hodge filtration $\F$ of $\V$ to $V = \V_{s_o}$ via the
Gauss--Manin connection $\nabla$ of $\V$.
\vskip 3pt

\par As in the pure case, the classifying spaces $\M$ defined above are
complex manifolds upon which a real Lie group acts transitively by complex
automorphisms.  In this subsections below, we shall introduce a certain
\lq\lq maximally homogeneous\rq\rq{} hermitian metric on $\M$, and compute
its curvature.

\proclaim{Theorem 2.3 [P2]} The classifying space $\M$ is a complex manifold 
upon which the real Lie group
$$
	\G = \{\, g\in GL(V)^W \mid Gr(g)\in Aut_{\R}(Q) \,\}
$$
acts transitively by automorphisms, where $GL(V)^W$ denotes the stabilizer of 
$W$ in $GL(V)$, and $Gr(g)$ denotes the induced action of $g\in GL(V)$ on 
$Gr^W$.
\endproclaim
\demo{Proof} That $\G$ acts transitively on $\M$ is a matter of simple linear 
algebra. In particular, since $\G$ acts transitively on $\M$, the orbit 
$\check\M\subseteq X$
of $F_o\in\M$ under the action of the complex Lie group
$$
	\Gc = \{\, g\in GL(V)^W \mid Gr(g)\in Aut_{\C}(Q) \,\}
$$
is well defined, independent of $F_o$.  Therefore, in order to show that $\M$ 
is a complex manifold on which $\G$ acts by automorphisms, it is sufficient 
to show (cf\. [P2]) that $\M$ is an open subset of 
$\check\M\cong\Gc/\Gc^{F_o}$, i.e\. for every $F\in\M$, there exists a 
neighborhood $U$ of $1$ in $\Gc$ such that 
$$
	g_{\C}\in U \implies g_{\C}.F\in\M
$$
\enddemo

\par In order to construct a hermitian metric on $\M$, we now recall the
following result from [CKS]:

\proclaim{Theorem 2.4} Let $(F,W)$ be a mixed Hodge structure.  Then, 
there exists a unique, functorial bigrading
$$
	V = \bigoplus_{p,q}\, I^{p,q}				\tag{2.5}
$$
of the underlying vector space $V=V_{\R}\otimes\C$ such that
\roster
\item"(a)" $F^p = \oplus_{a\geq p}\, I^{a,b}$;
\item"(b)" $W_k = \oplus_{a+b\leq k}\, I^{a,b}$;
\item"(c)" $\overline{I^{p,q}} = I^{q,p} \mod \oplus_{r<q,s<p}\, I^{r,s}$.
\endroster
\endproclaim

\proclaim{Corollary 2.6} Each choice of graded-polarization $Q = \{Q_k\}$ of 
$(F,W)$ determines a unique, functorial mixed Hodge metric $h_F$ on the 
underlying vector space $V$ such that
\roster
\item"(i)"  The decomposition $(2.5)$ is orthogonal with respect to $h_F$;
\item"(ii)" $u$, $v\in I^{p,q}\implies h_F(u,v) 
		= i^{p-q}Q_{p+q}([u],[\bar v])$.
\endroster
Accordingly, via the standard identification of $T_F(\M)$ with
a subspace of 
$$
	T_F(X) = \bigoplus_p\, \Hom(F^p,V/F^p)
$$
the mixed Hodge metric $(2.6)$ extends to a hermitian metric $h$ on $T(\M)$.
\endproclaim

\remark{Remark} Equivalently, the induced metric $(2.6)$ on $T(\M)$ can be
described as follows:  Let $F$ be a point in $\M$.  Then, application of 
Theorem $(2.4)$ to the mixed Hodge structure $(F^.\gg_{\C},W_.\gg_{\C})$ 
defines a functorial bigrading 
$$
	\gg_{\C} = \bigoplus_{r+s\leq 0}\, \gg^{r,s}_{(F,W)}	\tag{2.7}
$$
such that
$$
	{\frak t}_F = \bigoplus_{r<0}\,\gg^{r,s}_{(F,W)}	
$$ 
is a vector space complement to the isotopy algebra $\gg_{\C}^F$
of $F$ in $\gg_{\C}$.  Consequently,
$$
	T_F(\M) \cong {\frak t}_F				\tag{2.8}
$$
via the differential of the exponential map
$$
	e:{\frak t}_F\to\check\M,\qquad e(u) = \exp(u).F
$$
Moreover, relative to the isomorphism $(2.8)$,
$
	h_F(\a,\b) = Tr(\a\b^*)		
$.
\endremark
\vskip 3pt

\par In the pure case, the metric $(2.6)$ can be identified with a 
$\G$-invariant metric on the corresponding classifying space of pure, 
polarized Hodge structure $\Cal D$.  In contrast, in the mixed case, 
the action of $\G$ on $\M$ usually has 
non-compact isotopy, and hence there usually do not exist any $G$-invariant 
metrics on $\M$.   Nonetheless, both the decomposition $(2.5)$ and the metric 
$(2.6)$ are maximally homogeneous in the following sense:

\proclaim{Theorem 2.9 [K]} Let $F\in\M$, $G_{\R} = \G\cap GL(V_{\R})$ and 
$$
	\lam_{(F,W)} = \bigoplus_{r,s<0}\, \gg^{r,s}_{(F,W)}
$$
Then, 
$$
	\M = \G_{\R}\exp(\lam_{(F,W)}).F			\tag{2.10}
$$
Moreover, given any element $g\in\G_{\R}\cup \exp(\lam_{(F,W)})$:
\roster
\item"(i)"  $I^{p,q}_{(g.F,W)} = g.I^{p,q}_{(F,W)}$;
\item"(ii)" The induced map $L_{g*}:T_F(\M)\to T_{g.F}(\M)$ is an isometry.
\endroster
\endproclaim

\par To compute the curvature of $T(\M)$ with respect to the
mixed Hodge metric, let us fix a point $F\in\M$.  Then, on account of equation
$(2.10)$, every element $g_{\C}\in\Gc$ such that $g_{\C}.F\in\M$ admits a 
factorization of the form:
$$
	g_{\C} = g_{\R}e^{\l}f					\tag{2.11}
$$
where $g\in\G_{\R}$, $e^{\l}\in\exp(\lam_{(F,W)})$ and $f\in\Gc^F$.  Moreover,
(cf\. [P2]) by restricting the possible values of $\l$ and $\log(f)$ one can
define a distinguished real-analytic factorization of the form 
$(2.11)$ over a neighborhood of $1\in\Gc$.  Accordingly, by combining this
factorization with Theorem $(2.9)$, we can then calculate the curvature
of $\M$ following [D1]:

\proclaim{Theorem 2.12 [P1]} Let $F\in\M$, and 
$
	\gg_{\C} = \np\oplus\nz\oplus\nn\oplus\lam		
$
denote the decomposition of $\gg_{\C}$ defined by the subalgebras
$$
\aligned
	\np &= \bigoplus_{r\geq 0,\, s<0}\, \gg^{r,s}_{(F,W)}	\\
	\nz &= \gg^{0,0}_{(F,W)}
\endaligned\qquad
\aligned
	\nn  &= \bigoplus_{r<0,\, s\geq 0}\, \gg^{r,s}_{(F,W)}	\\
	\lam &= \bigoplus_{r,s<0}\, \gg^{\r,s}_{(F,W)}
\endaligned							
$$
Let $\pip$, $\piz$, $\pin$ and $\pil$ denote the corresponding projection
operators form $\gg_{\C}$ onto $\np$, $\nz$, $\nn$ and $\lam$.  Then, 
relative to the identification $(2.8)$, the hermitian holomorphic 
curvature of $T(\M)$ at $F$ with respect to the mixed Hodge metric $(2.6)$
is given by the formula:
$$
	R(u,v) = S(u,\bar v) - S(v,\bar u)		
$$
where
$$
\aligned
	S(u,\bar v) 
	&= \piq\,\ad((\pip[\bar v,u] + \half\piz[\bar v,u])
		 +(\pip[\bar u,v] + \half\piz[\bar u,v])^*)		\\
	&\hph{aa}+[\piq\,\ad \pip(\bar v),\piq\,\ad \pip(\bar u)^*]
\endaligned
$$
and $\piq$ denotes orthogonal projection from $gl(V)$ onto ${\frak t}_F$
with respect to $h_F$.
\endproclaim

\proclaim{Corollary 2.13} The holomorphic sectional curvature of $\M$
along $u\in T_F(\M)$ is given by the formula
$
	R(u) = h_F(S(u,\bar u)u,u)/h_F^2(u,u)
$.
\endproclaim

\remark{Remark} Unlike the pure case, the mixed Hodge metric $h$ need not
have negative holomorphic sectional curvature along horizontal directions.
The underlying reason for this is that $\G$ need not be semisimple, and
hence one can construct holomorphic, horizontal maps $F:\C\to\M$.
\endremark
\vskip 3pt

\par Following [K], in order to address the fact that $\G$ usually acts
with non-compact isotopy on $\M$, we now construct a natural fibration 
$\M\to\M_{\R}$ such that:
\roster
\item"(i)" $\G_{\R}$ acts transitively by isometries on $M_{\R}$;
\item"(ii)" The fiber over $\hat F$ is isomorphic to the subalgebra
$$
	\lam_{(\hat F,W)}\cap Lie(\G_{\R})
$$ via the map $\l \mapsto e^{i\l}.\hat F$.
\endroster

\par To this end, we recall that a grading of an increasing filtration $W$ of 
a finite dimensional vector space $V$ is a semisimple endomorphism $Y$ of $V$ 
such that $W_k$ is the direct sum of $W_{k-1}$ and the $k$-eigenspace $E_k(Y)$
for each index $k$.  In particular, by Theorem $(2.4)$, each mixed Hodge 
structure $(F,W)$ induces a functorial grading $Y=Y_{(F,W)}$ on the underlying
weight filtration $W$ via the rule:
$$
	E_k(Y) = \bigoplus_{p+q=k}\, I^{p,q}			\tag{2.14}
$$
Accordingly, a mixed Hodge structure $(F,W)$ is said to be {\it split over}
$\R$ if and only if the associated grading $(2.14)$ is defined over $\R$,
i.e\. $\overline{I^{p,q}} = I^{q,p}$.

\proclaim{Theorem 2.15 [K]} The locus of points $F\in\M$ such that $(F,W)$
is split over $\R$ is a $C^{\infty}$ submanifold of $\M$ on which $\G_{\R}$
acts transitively by isometries.
\endproclaim

\par To continue [K], let $\pi:\M\to\M_{\R}$ be a $C^{\infty}$ fibration such 
that:
\roster
\item"(a)" $\pi(F)\in\exp(\lam_{F,W}).F$; 
\item"(b)" $g\in\G_{\R}\implies\pi(g.F) = g.\pi(F)$;
\item"(c)" $F\in\M_{\R}\implies \pi(F) = F$.
\endroster
Then, on account of the fact that
$$
	\exp(\lam_{(F,W)})\cap\G^F = 1,
$$
the equation 
$$
	\pi(F) = e(F)^{-1}.F
$$
defines a $C^{\infty}$ function $e:\M\to\G$ 
such that
\roster
\item $e(F)\in\exp(\lam_{(F,W)})$;
\item $\hat F := e(F)^{-1}.F \in\M_{\R}$;
\item $g\in\G_{\R} \implies e(g.F) = \Ad(g)e(F)$;
\item $F\in\M_{\R}\implies e(F) = 1$.
\endroster
Conversely, given a $C^{\infty}$ function $e:\M\to\G$ which satisfies
conditions $(1)$--$(4)$, the above process can be inverted to define a 
corresponding fibration $\pi:\M\to\M_{\R}$ as above.  Thus, as a consequence 
of the next result, there exists a unique fibration $\M\to\M_{\R}$ such that
$$
	\overline{e(F)} = e(F)^{-1}
$$

\proclaim{Theorem 2.16 [CKS]} Let $(F,W)$ be a mixed Hodge structure.  Then,
there exists a unique, real element 
$$
	\delta\in \lam_{(F,W)} = \bigoplus_{r,s<0}\, gl(V)^{r,s}_{(F,W)}
$$
such that $(\hat F,W) = (e^{-i\delta}.F,W)$ is split over $\R$.
\endproclaim
\demo{Proof} Let $Y =Y_{(F,W)}$ denote the grading $(2.14)$ of $W$.  Then,
by virtue of Theorem $(2.4)$, 
$$
	\bar Y = Y \mod \lam_{(F,W)}
$$
Consequently (cf\. [CKS]), there exists a unique real element $\delta$ of 
$\lam_{(F,W)}$ such that
$$
	\bar Y = e^{-2i\delta}.Y
$$
Therefore, by virtue of part $(\text{i})$ of Theorem $(2.9)$, 
$(\hat F, W) = (e^{-i\delta}.F,W)$ is split over $\R$, with grading 
$
	Y_{(\hat F,W)} = e^{-i\delta}.Y_{(F,W)}
$
\enddemo

\par In particular, since both the mixed Hodge metric and the splitting
operation $(2.16)$ depend upon the Deligne--Hodge decomposition $(2.5)$,
the complexity of the fibration $\M\to\M_{\R}$ provides a measure of 
the failure of $\G$ to act on $\M$ by isometries.  As such, the next result 
implies that the geometry of the classifying spaces considered in \S 1 
should be \lq\lq simple\rq\rq{} [cf\. Theorem $(2.19)$]:

\proclaim{Theorem 2.17} Let $\M$ be a classifying space of type 
$(\text{\rm I})$ or $(\text{\rm II})$ \cite{cf\. \S 1}, and 
$$
	Lie_{-r}(W) = \{\, \a\in gl(V) \mid \a(W_k)\subseteq W_{k-r}\,\}
$$
Then, the fibration $\M\to\M_{\R}$ defined by Theorem $(2.16)$ is isomorphic
to the trivial fibration
$$
	\M \cong\R^d\times\M_{\R}
$$
where $d = \dim_{\C}\,Lie_{-2}(W)$.
\endproclaim
\demo{Proof} If $\M$ is type $(\text{I})$ then $d=0$ and every point $F\in\M$ 
is split over $\R$ due to the short length of $W$.
Similarly, if $\M$ is type $(\text{II})$ then
$$
	\lam_{(F,W)} = \gg^{-1,-1} = \bigoplus_{p+q=-2}\,\gg^{r,s} 
		     = Lie_{-2}(W)					
		     						\tag{2.18}
$$
due to the Hodge numbers of $\M$.  Consequently, in this case, the 
fibration $\M\to\M_{\R}$ is given by the formula
$$
	e^{i\l}.F \mapsto F,\qquad 
	\l\in Lie_{-2}(W)\cap gl(V_{\R}),\quad F\in\M_{\R}	
$$
\enddemo

\proclaim{Theorem 2.19} Let $\M$ be a classifying space of type 
$(\text{\rm I})$ or $(\text{\rm II})$.  Then, the subgroup 
$$
	\tlG = \{\, g\in\G \mid Gr(g)\in Aut_{\R}(W_k/W_{k-2})\,\}
$$
of $\G$ consisting of those elements $g\in\G$ which induce real automorphisms 
of $W_k/W_{k-2}$ for all $k$, acts transitively on $\M$ by isometries.
\endproclaim
\demo{Proof} If $\M$ is type $(\text{I})$ then $\M=\M_{\R}$ and 
$\tlG=\G_{\R}$, so we're done by Theorem $(2.15)$.  Suppose therefore that 
$\M$ is type $(\text{II})$.  Then, since $\tlG$ contains the subgroups 
$\G_{\R}$ and 
$$
	\exp(\lam_{(F,W)}) = \exp(Lie_{-2}(W))
$$
for every point $F\in\M$, it then follows from Theorem $(2.9)$ that $\tlG$
acts transitively on $\M$.  To see that $\tlG$ acts by isometries, recall
[CKS] that the set ${\Cal Y}(W)$ consisting of all gradings $Y$ of $W$
is an affine space upon $\exp(Lie_{-1}(W))$ acts simply
transitively by the rule
$$
	g.Y = \Ad(g)Y						\tag{2.20}
$$
Accordingly, given any element $g\in\G$ and any grading $Y\in\Cal Y(W)$,
there exists unique elements $g^Y\in\G^Y$ and $g_{-1}\in\exp(Lie_{-1}(W))$
such that 
$$
	g = g_{-1}g^Y						\tag{2.21}
$$
and $g.Y = g_{-1}.Y$.

\par Suppose now that $Y=\bar Y$.  Then, since every element of $\G$ acts
by real automorphisms on $Gr^W$, the corresponding factor $g^Y$ appearing
in $(2.21)$ actually belongs to $\G_{\R}$.  Furthermore, since $\M$ is type 
$(\text{II})$, $g_{-1}=e^{\a}$ can be factored as
$$
	g_{-1} = (1 + \a_{-1})(1 + \a_{-2})			\tag{2.22}
$$
where $\a_{-j}\in E_{-j}(\ad Y)$.  In particular, if $g\in\tlG$ then 
$\a_{-1}\in gl(V_{\R})$ since $g=g_{-1}g^Y$ acts by real automorphisms 
on $W_k/W_{k-2}$.  Consequently,
$$
	g = g_{-1}g^Y = \{(1+\a_{-1})g^Y\}\{(g^Y)^{-1}(1+\a_{-2})g^Y\}
								\tag{2.23}
$$
where the first term in curly braces on the right hand side of $(2.23)$
belongs to $\G_{\R}$, while the second term belongs $\exp(Lie_{-2}(W))$.
Therefore, by Theorem $(2.9)$ and equation $(2.18)$, 
$$
	L_{g*}:T_F(\M)\to T_{g.F}(\M)
$$
is an isometry for all $F\in\M$.
\enddemo

\remark{Remark} The proof of Theorem $(2.19)$ implies the following
additional fact:  If $\M$ is type $(\text{I})$ or $(\text{II})$ then
$h\in\tlG$, $F\in\M\implies I^{p,q}_{(h.F,W)} = h.I^{p,q}_{(F,W)}$.
\endremark

\head \S 3.\quad Limits of Mixed Hodge Structure \endhead

\par Let $\V\to\Delta^*$ be a variation of graded-polarized mixed
Hodge structure.  Then, in contrast to the pure case, the period map of $\V$
can have irregular singularities at the origin.  The source of this apparent 
disparity lies in the geometry of  the associated classifying spaces.  
Namely, unlike the pure case [S], the classifying spaces of graded-polarized 
mixed Hodge structure $\M$ discussed in \S 2 need not have negative 
holomorphic sectional curvature along horizontal directions.

\par Nevertheless, by comparison with the $\ell$-adic case, Deligne
conjectured in [D3] that the period map of a variation of mixed Hodge
structure arising from a family of complex algebraic varieties should
not have such irregular singularities.  Furthermore, according to [D3],
there should exist a category of \lq\lq good\rq\rq{} variations of 
mixed Hodge structures which both contains all of the geometric
variations and possesses the following salient features of the pure case:
\roster
\item"(a)" The existence of the limiting mixed Hodge structure;
\item"(b)" In the geometric case, the limiting Hodge structure $(a)$ should
admit a de Rham theoretic construction in terms of the log complex of the 
underlying morphism $f:X\to\Delta$;
\item"(c)" The existence of a functorial mixed Hodge structure on the 
cohomology $H^*(X,\V)$ of a good variation $\V\to X$;
\item"(d)" Nilpotent Orbit Theorem: The period map of a good variation
of mixed Hodge structure should be asymptotic to the corresponding nilpotent
orbit.
\endroster

\par In [SZ], Steenbrink and Zucker formulated the following definition of
a good variation:

\definition{Definition 3.1} A variation of graded-polarized mixed Hodge
structure $\V\to\Delta^*$ with unipotent monodromy is admissible
if  
\roster
\item"(i)"  The limiting Hodge structure $F_{\infty}$ of $\V$ exists;
\item"(ii)" The relative weight filtration $\rel W = \rel W(N,W)$ exists.
\endroster
\enddefinition 


\par The first evidence that this is indeed the correct definition is 
Deligne's proof in the appendix to [SZ] that conditions $(\text{i})$ and 
$(\text{ii})$ already imply that $(a)$ the pair $(F_{\infty},\rel W)$ is a 
mixed Hodge structure, relative to which $N$ is a $(-1,-1)$-morphism.  
Additional evidence is provided by the following two results [E] [Sa2], 
special cases of which are proven in [SZ]:
\roster
\item"---" Every geometric variation is admissible, and admits a de Rham
theoretic construction $(b)$ of its limiting mixed Hodge structure 
$(F_{\infty},\rel W)$;
\item"---" The cohomology $H^*(X,\V)$ of an admissible variation $\V\to X$
admits a functorial mixed Hodge structure $(c)$.
\endroster


\par In this section, we consider the singularities $(d)$ of the period map
$$
	\varphi:\Delta^*\to\Gamma\backslash\M			\tag{3.2}
$$
of an admissible variation $\V\to\Delta^*$ with unipotent monodromy.  To 
this end, let $p:U\to\Delta^*$ denote the universal cover of the punctured
disk by the upper half-plane, and $(s,z)$ be a pair of coordinates relative
to which $p$ assumes the form $s=e^{2\pi iz}$.  Then, by virtue of the
local liftablity of $\varphi$, there exists a holomorphic, horizontal map
$F:U\to\M$ which makes the following diagram commute:
$$
\CD
	U @> F >> \M						\\
	@V p VV	@VVV						\\
	\Delta^* @> \varphi >> \Gamma\backslash \M
\endCD\tag{3.3}
$$ 
Consequently, by the commutativity of $(3.3)$, the function 
$$
	\psi(z) := e^{-zN}.F(z)					\tag{3.4}
$$
descends to a well defined map $\psi(s):\Delta^*\to\check\M$.  Moreover,
we have the following result:

\proclaim{Lemma 3.5} $\V$ is admissible if and only if both the relative
weight filtration $\rel W$ and the limiting Hodge filtration
$$
	F_{\infty} = \lim_{s\to 0}\,\psi(s)
$$
exist.
\endproclaim

\par Thus, by the theorem of Deligne [SZ] quoted above, given an admissible 
variation $\V\to\Delta^*$, each choice of coordinates $(s,z)$ as above 
defines an associated limiting mixed Hodge structure $(F_{\infty},\rel W)$.  
Furthermore, just as in \S 2, the pair $(F_{\infty},\rel W)$ induces a 
functorial decomposition 
$$
	\gg_{\C} = \bigoplus_{r,s}\, \gg^{r,s}_{(F_{\infty},\rel W)}
								\tag{3.6}
$$
such that
$$
	{\frak t}_{\infty} 
	= \bigoplus_{r<0}\, \gg^{r,s}_{(F_{\infty},\rel W)}
$$
is a vector space complement to the isotopy algebra $\gg_{\C}^{F_{\infty}}$
in $\gg_{\C}$.  As such, near $s=0$, 
$$
	\psi(s) = e^{\Gamma(s)}.F_{\infty}			
$$
relative to a unique ${\frak t}_{\infty}$-valued holomorphic function 
$\Gamma(s)$ such that $\Gamma(0) = 0$.  Accordingly, by the definition of 
$\psi(s)$, 
$$
	F(z) = e^{zN}e^{\Gamma(s)}.F_{\infty}			\tag{3.7}
$$
for $\text{Im}(z)>>0$.  Moreover, just as in the pure case the period map
$F(z)$ is asymptotic to the associated nilpotent orbit 
$
	\theta(z) = e^{zN}.F_{\infty}
$ 
obtained by setting $\Gamma(s)=0$ in equation $(3.7)$: 
\vskip 3pt

\definition{Definition 3.8} An admissible, 1-variable nilpotent orbit 
is a holomorphic map $\theta:\C\to\check\M$ of the form
$$
	\theta(z) = e^{zN}.F
$$
where $F\in\check\M$ and $N$ is a nilpotent element of $\gg_{\R}$ such 
that
\roster
\item"---" $N(F^p)\subseteq F^{p-1}$;
\item"---" $\theta(z)\in\M$ for $\text{Im}(z)>>0$;
\item"---" (Admissibility): The relative weight filtration $\rel W(N,W)$ 
exists.
\endroster
\enddefinition

\proclaim{Theorem 3.9 (Nilpotent Orbit Theorem) [P3]} Let $\V\to\Delta^*$
be an admissible variation of graded-polarized mixed Hodge structure with
unipotent monodromy.  Then,
\roster
\item $\theta(z) = e^{zN}.F_{\infty}$ is an admissible nilpotent orbit;
\item There exists non-negative constants $\a$, $\b$ and $K$ such that 
$\text{Im}(z)>\a\implies\theta(z)\in\M$ and 
$$
	d_{\M}(F(z),\theta(z)) < K \text{Im}(z)^{\b}e^{-2\pi\text{Im}(z)}
$$
\endroster
\endproclaim

\par The proof of Theorem $(3.9)$ depends upon the follow results 
[S] [CKS] [D4] about split orbits which play a fundamental role
in \S 4--9:

\definition{Definition 3.10} A split orbit is an admissible nilpotent orbit
$(e^{zN}.\hat F,W)$ for which the associated limiting mixed Hodge structure
$(\hat F,\rel W)$ is split over $\R$.
\enddefinition

\par In the pure case, the notion of split and $\text{\rm SL}_2$-orbit 
coincide.Therefore, by [S] [CKS] we have the following classification of such 
orbits:

\definition{Definition 3.11} Let $H$ be a pure Hodge structure of weight $k$,
and $e=(1,0)$ and $f=(0,1)$ denote the standard  basis of $\C^2$.  Define 
$S(1)$ to be the standard representation of $sl_2(\C)$ on $\C^2$ equipped 
with the pure Hodge structure of weight one obtained by declaring
$$
	\nu_+ = e + i f,\qquad \nu_- = e - if			\tag{3.12}
$$
to be of type $(1,0)$ and $(0,1)$ respectively.  Then, a representation
of $sl_2(\C)$  on $H$ is Hodge if it induces a morphism of Hodge
structures from $sl_2(\C)\subset S(1)\otimes S(1)^*$ to 
$End(H) = H\otimes H^*$. 
\enddefinition

\remark{Remark} In [S], $\text{\rm SL}_2(\R)$ acts on upper half-plane via the
rule $z\mapsto (az + b)/(cz +d)$.  In [CKS], $\text{\rm SL}_2(\R)$ acts by
$z\mapsto (c+dz)/(a+bz)$.  We follow [CKS].
\endremark
 
\proclaim{Theorem 3.13 [CK], [CKS]} Let $\Cal D$ be a classifying space of 
pure Hodge structure, $F_o\in\Cal D$ and $\psi:\text{\rm SL}_2(\R)\to\G_{\R}$ 
be a representation of $\text{\rm SL}_2(\R)$.  Then,
$$
	\theta(g.\sqrt{-1}) = \psi(g).F_o
$$
is an $\text{\rm SL}_2$-orbit if and only if $\rho = \psi_*$ is Hodge with 
respect to $F_o$.
\endproclaim

\proclaim{Theorem 3.14 [S]} Let $H$ be a Hodge representation and 
$S(k) = Sym^k(S(1))$.  Then, $H$ can be decomposed into a direct sum 
of irreducible Hodge submodules.  Furthermore, every irreducible Hodge 
representation is isomorphic to one of the following types\footnote{By
convention $S(0)=H(0)$.}
\roster
\item"(a)" $H(d)\otimes S(m)$, $m\geq 0$;
\item"(b)" $E(p,q)\otimes S(n)$, $p-q>0$, $n\geq 0$;
\endroster
where $H(d)=\C$ and $E(p,q)=\C^2$ denote the following Hodge structures,
equipped with the trivial action of $sl_2(\C)$:
\roster
\item"---" $H(d)$ is weight $-2d$ and type $(-d,-d)$;
\item"---" $E(p,q)$ is weight $p+q$, $\nu_+$ of type $(p,q)$
and $\nu_-$ of type $(q,p)$.
\endroster
\endproclaim

\remark{Remark} Let $H$ be a Hodge representation, and $Q$ be a polarization
of $H$ which is compatible with the given action of $sl_2(\C)$.  Then, the
decomposition of Theorem $(3.14)$ can be chosen to be orthogonal with
respect to $Q$.  Furthermore, each irreducible summand is isomorphic to 
one of the standard tensor products $(a)$ $(b)$ equipped with the following
polarizations: 
\roster
\item"---" $H(d):Q(1,1) = 1$;
\item"---" $S(1):Q(e,f) = 1$, $S(k) = Sym^k(S(1))$;
\item"---" $E(p,q):Q(e,f) = i^{q-p+1}$.
\endroster
\endremark
\vskip 3pt

\par In the mixed case, a split orbit $\theta(z) = e^{zN}.\hat F$ induces 
$\text{\rm SL}_2$-orbits on $Gr^W$.  Accordingly, each choice of grading $Y$ 
of $W$ defines a corresponding lift of the associated representations of 
$sl_2$ on $Gr^W$ to a representation
$$
	\rho_Y:sl_2(\C)\to\gg_{\C}
$$
In [D4], Deligne showed how to use the limiting mixed Hodge structure of
$\theta(z)$ to make a distinguished choice of grading $Y$ such that the
associated representation $\rho_Y$ has a number of very special properties.
To state Deligne's result, let
$$
	\n_o = \pmatrix 0 & 0 \\ 1 &  0 \endpmatrix,\qquad
	\h =  \pmatrix  1 & 0 \\ 0 & -1 \endpmatrix,\qquad
	\n_o^+ = \pmatrix 0 & 1 \\ 0 & 0 \endpmatrix		\tag{3.15}
$$
denote the standard generators of $sl_2(\C)$ and $\rel Y$ denote the 
grading $(2.14)$ of the relative weight filtration $\rel W$ defined 
by the $I^{p,q}$'s of the limiting mixed Hodge structure of $\theta$.

\proclaim{Theorem 3.16 [D4]} Let $\theta(z) = e^{zN}.\hat F$ be a split 
orbit.  Then, there exists a unique, functorial $\R$-grading $Y$ of $W$ such 
that
\roster
\item $[\rel Y,Y] = 0$;
\item $[N - \rho_Y(n_o),\rho_Y(n_o^+)] = 0$.
\endroster
Furthermore, if 
$$
	N  = N_0 + N_{-1} + N_{-2} + \cdots			\tag{3.17}
$$
denotes the decomposition of $N$ with respect to the eigenvalues of 
$\ad Y$ and 
$$
	N_0 = \rho(\n_o),\qquad H = \rho(\h),
	\qquad N_0^+ = \rho(\n_o^+)				\tag{3.18}
$$
denotes the $sl_2$-triple defined by the representation $\rho = \rho_Y$ then:
\roster
\item"(a)" $N_{-k}$, $k>0$ is highest weight $k-2$ with respect to $\rho$;
\item"(b)" $H = \rel Y - Y$;
\item"(c)" $e^{zN_0}.\hat F$ is an $\text{\rm SL}_2$-orbit ({\rm Data:} 
$F_o = e^{iN_0}.\hat F$, $\psi_* = \rho$);
\item"(d)" $Y$ preserves $\hat F$, $Y_{(e^{iyN_0}.\hat F,W)} = Y$, and
$Y_{(e^{zN}.\hat F,W)} = e^{zN}.Y$;
\endroster
In particular, as consequence of $(a)$, $N_{-1} = 0$ and $[N_0,N_{-2}] = 0$.
\endproclaim
\demo{Proof} See [KP] [P3] [Sch].
\enddemo

\remark{Remark} More generally, in [D4] Deligne proved the following 
result:  Let $\rel Y$ be a grading of the relative weight filtration
such that $[\rel Y,N] = -2N$.  Assume $\rel Y$ preserves $W$.  Then,
there exists a system of graded representations $Gr(\rho)$ and a unique
functorial $\C$-grading 
$$
	Y = Y(N,\rel Y)						\tag{3.19}
$$
of $W$ which satisfies conditions $(1)$--$(2)$ and $(a)$--$(b)$ of 
Theorem $(3.16)$.  Accordingly, if $(e^{zN}.F,W)$ is an admissible
nilpotent orbit then application of $(3.19)$ to $N$ and 
$\rel Y = Y_{(F,\rel W)}$ defines a corresponding grading
$$
	Y = Y(F,W,N)						\tag{3.20}
$$
of $W$.
\endremark
\vskip 3pt

%
%


\comment

\endcomment




\head \S 4.\quad $\text{\rm SL}_2$-Orbit Theorem \endhead

\par Let $X$ be a complex algebraic variety.  Then, by [D2, III, \S 8.2] the 
hodge numbers $h^{p,q}$ of the mixed Hodge structure attached to $H^n(X,\C)$ 
satisfy the following numerical conditions:
\roster
\item"(i)" $h^{p,q} = 0$ unless $0\leq p, q\leq n$;
\item"(ii)" If $X$ is proper, then $h^{p,q} = 0$ unless $p+q\leq n$;
\item"(iii)" If $X$ is smooth, then $h^{p,q} = 0$ unless $p+q\geq n$;
\item"(iv)" If $N = \dim(X)$ and $n\geq N$, then $h^{p,q} = 0$ unless
           $n-N\leq p,q\leq N$.
\endroster

\par Accordingly, by conditions $(\text{i})$ and $(\text{iv})$, 
given any complex algebraic variety $X$, the mixed Hodge structures attached 
to $H^1(X;\Z(1))$ and $H^{2N-1}(X;\Z(N))$ are of the form 
$$
	H_{\C} = I^{0,0}\oplus I^{0,-1}\oplus I^{-1,0}\oplus I^{-1,-1}
								\tag{4.1}
$$
with $Gr^W_{-1}$ polarizable, and hence determine [D2, III, \S 10.1] a
corresponding pair of 1-motives, called the Picard and Albanese 1-motives of 
$X$.  Likewise, given a family $f:X\to S$ of complex algebraic varieties, the 
local systems $Pic = R^1_{f*}(\Z(1))$ and $Alb = R^{2n-1}_{f*}(\Z(n))$ support
admissible variations of 1-motives of type $(\text{II})$ over a Zariski 
open subset of $S$.  Moreover, by conditions $(\text{ii})$ and 
$(\text{iii})$, $Pic$ and $Alb$ reduce to variations of type 
$(\text{I})$ whenever the generic fiber of $f$ is either proper or 
smooth.
\vskip 3pt

\par Returning now to the context of abstract variations, our main result
can be stated as follows:

\proclaim{Theorem 4.2 ($\text{\bf SL}_2$-Orbit Theorem)} Let $e^{zN}.F$ be an 
admissible nilpotent orbit of type $(\text{\rm I})$ or $(\text{\rm II})$, 
with relative weight filtration $\rel W = \rel W(N,W)$ and $\delta$-splitting
\cite{cf\. Theorem $(2.16)$}
$$
            (F,\rel W) = (e^{i\delta}.\hat F,\rel W)               
$$
Define $\hh = Lie(\tlG)$ \cite{cf\. Theorem $(2.19)$}.  Then, there exists an 
element 
$$
    \zeta\in\hh\cap\ker(N)\cap\lam_{(\hat F,\rel W)}
$$
and distinguished real analytic function $g:(a,\infty)\to\tlG$ such that 
\roster
\item"(a)" $e^{iyN}.F = g(y)e^{iyN}.\hat F$;

\item"(b)" $g(y)$ and $g^{-1}(y)$ have convergent series expansions about
$\infty$ of the form
$$
\aligned
       g(y) &= e^{\zeta}(1 + g_1 y^{-1} + g_2 y^{-2} + \cdots)      \\
  g^{-1}(y) &= (1 + f_1 y^{-1} + f_2 y^{-2} + \cdots)e^{-\zeta}         
\endaligned
$$
with $g_k$, $f_k\in \ker(\ad N_0)^{k+1}\cap\ker(\ad N_{-2})$;
\item"(c)" $\delta$, $\zeta$ and the coefficients $g_k$ are related by
the formula
$$
      e^{i\delta} 
      = e^{\zeta}\left(1 + \sum_{k>0}\, \frac{1}{k!}(-i)^k(\ad\, N_0)^k\,g_k
                 \right)
$$
\endroster
where $(N_0,H,N_0^+)$ denotes the $sl_2$ triple attached to $e^{zN}.\hat F$ 
by Theorem $(3.16)$, and $N = N_0 + N_{-2}$ is the corresponding decomposition
of $N$.  Moreover, $\zeta$ can expressed as a universal Lie polynomial over 
$\Q(\sqrt{-1})$ in the Hodge components $\delta^{r,s}$ of $\delta$ with 
respect to $(\hat F,W)$.  Likewise, the coefficients $g_k$ and $f_k$  
can be expressed as universal, non-commuting 
polynomials over $\Q(\sqrt{-1})$ in $\delta^{r,s}$ and $\ad\, N_0^+$.
\endproclaim

\par By way of applications of this result, we now state three general
consequences of Theorem $(4.2)$.  To this end, we note that, in conjunction 
with the nilpotent orbit theorem discussed in \S 3, one expects to be able to
reduce many questions regarding the asymptotic behavior of an admissible 
variation $\V\to\Delta^*$ to the case of split orbits via Theorem $(4.2)$.  
More precisely, one has:

\proclaim{Corollary 4.3} Let $\V\to\Delta^*$ be an admissible variation of 
type $(\text{\rm I})$ or $(\text{\rm II})$, with period map $F(z):U\to\M$ 
and nilpotent orbit $e^{zN}.F$.  Then, adopting the notation of Theorem 
$(4.2)$, there exists a distinguished, real--analytic function $\g(z)$ with 
values in $\hh$ such that, for $\Im(z)$ sufficiently large,
\roster
\item"(i)"  $F(z) = e^{xN}g(y)e^{iyN_{-2}}y^{-H/2}e^{\g(z)}.F_o$;
\item"(ii)" $|\g(z)| = O(\Im(z)^{\b}e^{-2\pi\Im(z)})$ as $y\to\infty$ and
$x$ restricted to a finite subinterval of $\R$, for some constant $\b\in\R$.
\endroster
where $F_o = e^{iN_0}.\hat F$.
\endproclaim
\demo{Proof} By equation $(3.7)$, we can write
$$
      F(z) = e^{zN}e^{\Gamma(s)}.F_{\infty},\qquad s = e^{2\pi i z}
$$
relative to a distinguished $\gg_{\C}$--valued holomorphic function 
$\Gamma(s)$ which vanishes at $s=0$.  
Therefore,
$$
\aligned
       F(z) &= e^{zN}e^{\Gamma(s)}.F 
             = e^{xN}e^{iyN}e^{\Gamma(s)}.F                         \\
	    &= e^{xN}e^{iyN}e^{\Gamma(s)}e^{-iyN}e^{iyN}.F
	     = e^{xN}e^{\Gamma_1(z)}e^{iyN}.F        
\endaligned                                                         
$$
where $\Gamma_1(z) = e^{iyN}e^{\Gamma(s)}e^{-iyN}$.  By Theorem $(4.2)$,
$$
       e^{iyN}.F = g(y)e^{iyN}.\hat F = g(y)e^{iyN_{-2}}y^{-H/2}.F_o
$$
since $y^{-H/2}.F_o = e^{iyN_0}.\hat F$.  Consequently, 
if  $h(y) = g(y)e^{iyN_{-2}}y^{-H/2}$ then 
$$
\aligned
        F(z) &= e^{xN}e^{\Gamma_1(z)}e^{iyN}.F 
              = e^{xN}e^{\Gamma_1(z)}h(y).F_o                       \\
             &= e^{xN}h(y)h^{-1}(y)e^{\Gamma_1(z)}h(y).F_o           
              = e^{xN}h(y)e^{\Gamma_2(z)}.F_o
\endaligned                                                         \tag{4.4}
$$
where $\Gamma_2(z) = h^{-1}(y)e^{\Gamma_1(z)}h(y)$.  Also, 
$$
         |\Gamma_2(z)| = O(\Im(z)^{\b}e^{-2\pi\Im(z)})              \tag{4.5}
$$
since $\Gamma(s)$ is a holomorphic function such that $\Gamma(0) = 0$, 
$e^{iyN}$ and $e^{iyN_{-2}}$ are polynomial in $y$, $g(y) = O(1)$ 
and $y^{H/2}$ acts as multiplication by an integral power of $y^{1/2}$ on
the eigenspaces of $H$.

\par To complete the proof, we now recall that by equation $(2.11)$,
we may write
$$
        e^{\Gamma_2(z)} = g_{\R}(z)e^{\l(z)}f(z)                    \tag{4.6}
$$
where each factor is real--analytic, and 
$$
      g_R(z)\in\G_{\R},\qquad \l(z)\in\lam_{(F_o,W)},\qquad f(z)\in\Gc^{F_o}
$$  
Accordingly, for $\Im(z)$ sufficiently large, there exists a unique 
$\hh$-valued function $\g(z)$ such that
$$
      e^{\g(z)} = g_{\R}(z)e^{\l(z)}
$$
By equation $(4.4)$, $\g(z)$ satisfies $(\text{i})$ since
$f(z)$ takes values in $\Gc^{F_o}$.  Likewise, $\g(z)$ satisfies condition
$\text{(ii)}$ by virtue of equation $(4.5)$ and the fact that the 
decomposition $(4.6)$ is real--analytic.
\enddemo

\remark{Remark} For variations of type $(\text{I})$, $N = N_0$.  
For variations of type $(\text{II})$,
$N=N_0 + N_{-2}$ and $\ker(N) = \ker(N_0)\cap\ker(N_{-2})$.
\endremark
\vskip 3pt

\par Our first application of Theorem $(4.2)$ is the following analog of the
1-variable norm estimates [S, Theorem $(6.6)$]:

\proclaim{Theorem 4.7} Let $\V\to\Delta^*$ be an admissible variation of 
type $(\text{\rm I})$ or $(\text{\rm II})$ with weight filtration $\W$ 
and relative weight filtration $\rel\W$.  Then, adopting the notation of 
Theorem $(4.2)$,
\roster
\item"(a)" The norm $||\sigma(s)||$ of a flat, global section of $\V$ remains
bounded as $s\to 0$;

\item"(b)" Over any angular sector $A$ of $\Delta^*$, a flat section $\sigma$
of $\rel\W_k$ satisfies the estimate
$$
                || \sigma(s) || = O((-\log|s|)^{\frac{k}{2}})
$$
provided $\W_{\ell} = 0$ for $\ell<0$.
\endroster
More generally, if $F(z):U\to\M$ denotes the period map of $\V$ then,
for $x=\text{Re}(z)$ restricted to a finite subinterval of $\R$, 
$$
    v\in E_k(H)\cap\ker(N_{-2}) \implies ||v||_{F(z)} = O(y^{\frac{k}{2}})    
								  \tag{4.8}
$$
as $y\to\infty$.
\endproclaim
\demo{Proof} The estimate $(4.8)$ implies items $(a)$ and $(b)$.  Indeed,
after pulling back $\V$ to the upper half-plane, a flat global section of 
$\V$ is represented by a constant vector\footnote{To show that 
$\ker(N) = \ker(N_0)\cap\ker(N_{-2})$ note that if $\V$ is type $(\text{I})$
then $N=N_0$, so there is nothing to prove.  If $\V$ is type $(\text{II})$
then $N_0$ must act trivially on $Gr^W_{2k}$ and $Gr^W_{2k-2}$ since they
are of pure Hodge type $(k,k)$ and $(k-1,k-1)$ respectively.} 
$$
	v\in\ker(N) = \ker(N_0)\cap\ker(N_{-2})			
$$  
Therefore, upon decomposing $v$ into its isotypical components with respect 
to the representation of $sl_2$ defined by $(N_0,H,N^+_0)$, it then follows 
that [since $N_{-2}$ commutes with $(N_0,H,N_0^+)$] each such component is 
also contained in $\ker(N_0)\cap\ker(N_{-2})$, and hence belongs
to $E_k(H)\cap\ker(N_{-2})$ for some index $k\leq 0$.  Consequently, by 
$(4.8)$, $||v||_{F(z)}$ is bounded.

\par Likewise, over any angular sector, a flat section of $\rel\W_k$ is
represented by a constant vector $v\in\rel W_k$.  Therefore, recalling
$(3.16b)$ that 
$$
          H = \rel Y - Y
$$
where $\rel Y$ is a grading of $\rel W$ and $Y$ is a grading of $W$
which commutes with $\rel Y$, it then follows that
$$
      W_{\ell} = 0\hph{a}\text{for}\hph{a} \ell<0 \implies
      \rel W_k \subseteq\bigoplus_{j\leq k}\, E_j(H)
$$
\comment
To verify this assertion, note that because $\rel Y$ and $Y$ commute, 
$$
       \rel W_k = \bigoplus_{a,b}\, \rel W_k\cap E_a(\rel Y)\cap E_b(Y)
$$
By assumption $E_b(Y) = 0$ if $b<0$. Likewise $\rel W_k\cap E_a(\rel Y) = 0$
for $a>k$ since $\rel Y$ is a grading of $\rel W$.  Consequently, 
$H = \rel Y - Y$ acts on any non-trivial summand 
$ \rel W_k\cap E_a(\rel Y)\cap E_b(Y)$ as multiplication by $a-b\leq k$.
\endcomment
Invoking $(4.8)$, one then obtains $(b)$.

\par To establish $(4.8)$, suppose that $\V$ is a split orbit, i.e.
$F(z) = e^{zN}.\hat F$.  Then, given a vector $v\in E_k(H)\cap\ker(N_{-2})$, 
$$
        || v ||_{e^{zN}.\hat F} 
      = || v ||_{e^{xN}e^{iyN}.\hat F}
      = || e^{-xN} v ||_{e^{iyN}.\hat F} 
      = || v + v'(x) ||_{e^{iyN}.\hat F}
$$
where 
$$
       v'(x) \in \bigoplus_{j\leq k-2}\, E_j(H)
$$
since $N_0:E_a(H)\to E_{a-2}(H)$, $N_{-2}(v)=0$, and  
$e^{xN} = e^{xN_0}e^{xN_{-2}}$ as $[N_0,N_{-2}] = 0$.
Accordingly, it suffices to show that 
$$
       v\in E_k(H)\cap\ker(N_{-2})\implies 
       ||v||_{e^{iyN}.\hat F} = y^{\frac{k}{2}} ||v||_{e^{iN}.\hat F}
$$
However, since $e^{zN}.\hat F$ is a split orbit, 
$$
       e^{iyN}.\hat F = e^{iyN_{-2}}y^{-H/2}e^{iN_0}.\hat F
$$
Therefore, as $H\in\gg_{\R}$ via $(3.14)$,
$N_{-2}\in\lam_{(\tilde F,W)}$ for all $\tilde F\in\M$ by $(2.18)$, and 
$v\in\ker(N_{-2})$,
$$
\aligned
      ||v||_{e^{iyN}.\hat F} 
      &= || v ||_{e^{iyN_{-2}}y^{-H/2}e^{iN_0}.\hat F}
       = || e^{-iyN_{-2}} v||_{y^{-H/2}e^{iN_0}.\hat F}      		\\
      &= || v||_{y^{-H/2}e^{iN_0}.\hat F} 
       = || y^{H/2} v||_{e^{iN_0}.\hat F} 
       = y^{\frac{k}{2}} || v ||_{e^{iN_0}.\hat F}
\endaligned
$$

\par More generally, given an admissible variation $\V\to\Delta^*$ of 
type $(\text{I})$ or $(\text{II})$, one can replicate the above 
argument mutatis mutandis using Corollary $(4.3)$.  The only trick is to 
note that since $f_k\in\ker(\ad\, N_0)^{k+1}$, the term 
$\Ad(y^{H/2})(f_k y^{-k})$ is at worst $O(1)$ in $y$, and 
$[N_{-2},g^{-1}(y)]=0$ since all the terms of the series expansion of 
$g^{-1}(y)$ belong to $\ker(\ad\, N_{-2})$.
\enddemo

\par Theorem $(4.7)$ shows that admissible variations of type 
$(\text{I})$ satisfy norm estimates which are identical to the
pure case.  The next result make a similar assertion regarding the 
holomorphic sectional curvature:

\proclaim{Theorem 4.9} Let $\V\to\Delta^*$ be an admissible variation of 
type $(\text{\rm I})$ with non-trivial monodromy logarithm $N$, and 
period map $F(z):U\to\M$.  Then, the holomorphic sectional curvature of 
$\M$ along $F(z)$ is negative, and bounded away from zero for $\Im(z)$ 
sufficiently large.
\endproclaim
\demo{Proof} By Corollary $(2.13)$, the holomorphic sectional curvature of 
$\M$ along $u\in T_F(\M)$ is given by a formula of the form
$$
        R(u) = \frac{h_F(S_F(u,\bar u)u,u)}{h^2_F(u,u)}
$$
relative to a $\G_{\R}$-invariant tensor field $S$.  Consequently, upon 
writing $F(z)$ 
$$
       F(z) = e^{xN}g(y)y^{-H/2}e^{\gamma(z)}.F_o
$$
as per Corollary $(4.3)$, one finds that [via the $\G_{\R}$-invariance of S]
$$
        R(F_*(d/dz)) 
	= \frac{h_{F_o}(S_{F_o}(\theta(z),\bar\theta(z))\theta(z),\theta(z))}
               {h_{F_o}(\theta(z),\theta(z))}                      \tag{4.10}
$$
where 
$$
        \theta(z) = \Ad(e^{-\g(z)})(\b^{-1,1}(y) + \b^{-1,0}(y))    \tag{4.11}
$$
and $\b^{-1,1}(y)$ and $\b^{-1,0}(y)$ denote the Hodge components of the
function 
$$
    \b(y) = \Ad(h^{-1}(y))N,\qquad h(y) = g(y)y^{-H/2}              \tag{4.12}
$$
with respect to the base point $F_o = e^{iN}.\hat F$.  In particular, as a 
consequence of the proof of Theorem $(4.2)$ for nilpotent orbits of type 
$(\text{I})$ given in \S 8,  $\b(y)$ admits a series expansion about infinity 
of the form
$$
    \b(y) = \sum_{n\geq 0}\, \b_n y^{-1-n/2}
$$
with leading order term $\b_0 = N$.  Therefore, by equations 
$(4.10)$--$(4.12)$, 
$$
      \lim_{\Im(z)\to\infty}\, R(F_*(d/dz)) 
      =  \frac{h_{F_o}(S_{F_o}(\xi,\bar\xi)\xi,\xi)}
              {h^2_{F_o}(\xi,\xi)}                                 \tag{4.13}
$$
where 
$$
        \xi = N^{-1,1} = \frac{1}{4}(iH + N_0 + N_0^+)             \tag{4.14}
$$
On the other hand, by Theorem $(2.12)$, 
$$
     h_{F_o}(S_{F_o}(\xi,\bar\xi)\xi,\xi)
     = -h_{F_o}([\bar\xi,\xi],[\bar\xi,\xi])  
$$
Thus,
$$ 
          \lim_{\Im(z)\to\infty}\, R(F_*(d/dz)) < 0.
$$
\enddemo

\remark{Remark} Theorem $(4.9)$ is false for variations of type 
$(\text{II})$.  In particular, if $\V$ is Hodge--Tate then
$R(F_*(d/dz)) = 0$ for all $z$.
\endremark
\vskip 10pt

\par To put the next result in context, we recall that in the pure case, 
Schmid's nilpotent orbit theorem asserts the existence of the limiting Hodge 
filtration of a variation of pure polarized Hodge structure $\V\to\Delta^*$.  
In the mixed case, this existence of the limiting Hodge filtration is assumed.
Less clear in the mixed case however is how the corresponding grading 
$$
       \Y(s) = Y_{(\F(s),\W)}
$$
of $\W$ behaves as $s\to 0$.

\proclaim{Theorem 4.15} Let $\V\to\Delta^*$ be an admissible variation 
of type $(\text{\rm I})$ or $(\text{\rm II})$ with period map $F(z):U\to\M$.  
Then, the limiting grading
$$
      Y_{\infty} = \lim_{\Im(z)\to\infty}\, e^{-zN}.Y_{(F(z),W)} 
$$
exists, and coincides with the grading $Y(F_{\infty},W,N)$ defined by 
equation $(3.20)$.
\endproclaim
\demo{Proof} By Corollary $(4.3)$,
$$
\aligned
      F(z) &= e^{xN}g(y)e^{iyN_{-2}}y^{-H/2}e^{\g(z)}.F_o             \\
           &= e^{xN}g(y)e^{\g_1(z)}e^{iyN_{-2}}y^{-H/2}.F_o
            = e^{xN}g(y)e^{\g_1(z)}e^{iyN}.\hat F
\endaligned                                                      
$$
where 
$$
         \g_1(z) = \Ad(e^{iyN_{-2}}y^{-H/2})\g(z)
$$ 
is a $\hh$-valued function of order 
$\Im(z)^{\b}e^{-2\pi\Im(z)}$, and $F_o = e^{iN_0}.\hat F$.  
Consequently, if 
$
       Y = Y(\hat F,W,N)
$
then 
$$
       e^{-zN}.Y_{(F(z),W)} 
       = e^{-iyN}g(y)e^{\g_1(z)}.Y_{(e^{iyN}.\hat F,W)}
       = e^{-iyN}g(y)e^{\g_1(z)}e^{iyN}.Y                        \tag{4.16}
$$
since 
$
       Y_{(e^{iyN}.\hat F,W)} = e^{iyN}.Y
$
by Theorem $(3.16d)$. 
Setting 
$$
       \g_2(z) = \Ad(e^{-iyN})\g_1(z)                             \tag{4.17}
$$
it then follows from equations $(4.16)$ and $(4.17)$ that
$$
    \lim_{\Im(z)\to\infty} e^{-zN}.Y_{(F(z),W)} 
    = \lim_{\Im(z)\to\infty} e^{-iyN}g(y)e^{iyN}e^{\g_2(z)}.Y     \tag{4.18}
$$
Therefore, by part $(b)$ of Theorem $(4.2)$,
$$
\aligned
    e^{-iyN}g(y)e^{iyN}
    &= e^{\zeta} e^{-iy\,\ad\, N} 
	   \left(1 + \sum_{k>0}\, g_k y^{-k} \right) \\ 
    &= e^{\zeta}\left(1 + \sum_{k>0}\sum_{j=}^k\, 
         \frac{1}{j!}(-i)^j (\ad N_0)^j\, g_k y^{j-k} \right)
\endaligned                                                     
$$
since $N = N_0 + N_{-2}$, $[N_0,N_{-2}] = 0$, 
$g_k\in\ker(\ad\, N_0)^{k+1}\cap\ker(\ad\, N_{-2})$ and 
$\zeta\in\ker(\ad\,N_0)\cap\ker(\ad\, N_{-2})$.
Consequently, by part $(c)$ of Theorem $(4.2)$,
$$
    \lim_{y\to\infty}\,e^{-iyN}g(y)e^{iyN}
    = e^{\zeta}\left(1 + \sum_{k>0}\, 
         \frac{1}{k!}(-i)^k (\ad N_0)^k\, g_k \right) = e^{i\delta}
	                                                       \tag{4.19}
$$
On the other hand, by equation $(4.17)$, $\g_2(z)$ is also of order
$\Im(z)^{\b}e^{-2\pi\Im(z)}$ for some constant $\b$.  Therefore, 
$$
     \lim_{\Im(z)\to\infty}\, e^{\g_2(z)} = 1                  \tag{4.20}
$$
Inserting equations $(4.19)$ and $(4.20)$ into equation $(4.18)$, it then
follows that 
$$
       Y_{\infty} = \lim_{\Im(z)\to\infty}\, e^{-zN}.Y_{(F(z),W)} 
        = e^{i\delta}.Y  = Y(F_{\infty},W,N)
$$
since 
$
    Y(F_{\infty},W,N) = e^{i\delta}.Y(\hat F_{\infty},W,N) 
    = e^{i\delta}.Y
$
by the functoriality of $Y$ (cf\. [P3]).
\enddemo

\remark{Remark} By [KP], Theorem $(4.15)$ is also true for unipotent 
variations (e.g. the variations attached to fundamental group of a smooth 
variety [HZ]) and variations for which the limiting mixed Hodge 
structure is split over $\R$ in some suitable coordinate system (e.g. the 
A--model variation considered in mirror symmetry [P2]).
\endremark

\head \S 5.\quad Arakelov Geometry \endhead

\par Let $M$ be a graded-polarized mixed Hodge structure.  Then, motivated
by the construction of [H] described below, we define the height of $M$
to be
$$
	h(M) = 2\pi||\delta||					\tag{5.1}
$$
where $\delta$ denotes the splitting of $M$ defined in \S 2, and $||*||$ 
denotes the mixed Hodge norm of $M$.  

\par To relate the height functional $(5.1)$ to the standard archimedean
height pairing defined by [Arak] [Beil] [GS], let $X$ be a non-singular
complex projective variety of dimension $n$, and $Z$ and $W$ be a pair
of algebraic cycles in $X$ of dimensions $d=\dim(Z)$ and $e=\dim(W)$
such that
\roster
\item"(i)"   $Z$ and $W$ are homologous to zero in $X$;
\item"(ii)"  $d + e = n-1$;
\item"(iii)" $|Z|\cap |W|=\emptyset$.
\endroster 
Then, as a consequence of \S 3 of [H], the mixed Hodge structure on 
$H_{2d+1}(X-|W|,|Z|;\Z(-d))$ carries a canonical subquotient $B=B_{Z,W}$
with graded pieces
$$
	Gr^W_0 \cong\Z(0),\quad 
	Gr^W_{-1}\cong H_{2d+1}(X;\Z(-d)),\quad
	Gr^W_{-2}\cong\Z(1)					\tag{5.2}
$$
such that
$$
	h(B_{Z,W}) = |\langle Z,W \rangle | 			\tag{5.3}
$$
where $\langle Z,W \rangle$ denotes the archimedean height of the pair 
$(Z,W)$.  

\par More precisely, via the cycles $Z$ and $W$, one obtains canonical
generators $1$ and $1'$ of $Gr^W_0(B)\cong\Z(0)$ and $Gr^W_{-2}(B)\cong\Z(1)$
respectively.  Moreover, as a consequence of Proposition $(3.2.13)$ in [H],
$$
	\delta(1) = \frac{1}{2\pi}\langle Z, W \rangle 1'	\tag{5.4}
$$
from which one then obtains equation $(5.3)$ via the definition of the
mixed Hodge metric.

\par Likewise, given a smooth, proper morphism $\pi:X\to S$ of relative 
dimension $n$, and a pair of flat, algebraic cycles $Z$ and $W$ in $X$ 
of relative dimensions $d$ and $e$ such that, for generic $s\in S$, $X_s$ 
is smooth and the triple $(X_s,Z_s,W_s)$ satisfies conditions 
$(\text{i})$--$(\text{iii})$ above, one obtains a corresponding height 
function
$$
	h(s) = \langle Z_s,W_s \rangle				\tag{5.5}
$$
over a Zariski dense open subset $S'$ of $S$.  

\par Let $D$ be a normal crossing divisor contained in the boundary
of a smooth partial compactification $\overline{S'}$ of $S'$.  In [H],
Hain analyzed the asymptotic behavior of $(5.5)$ near $D$ under the
assumption that the associated variation of mixed Hodge structure
$$
	\V\to S',\qquad \V_s = B_{Z_s,W_s}			\tag{5.6}
$$
induced constant variation of pure Hodge structure on $Gr^{\W}$.  In [L],
Lear computed the asymptotic behavior of $(5.5)$ under the assumption 
that $S$ is a curve using the theory of normal functions.  

\par In this section, we consider the asymptotic behavior of $(5.5)$ near
a normal crossing divisor $D$ about which $\V$ degenerates with unipotent
monodromy by applying Theorem $(4.2)$ to the 1-parameter 
degenerations $f^*(\V)$ obtained by pulling back $\V$ along a holomorphic 
map $f$ from the unit disk $\Delta$ into $\overline{S'}$.

\par To this end, let us assume for the moment that $\dim(S) = 1$ and $p$
is a point about which $\V$ degenerates with unipotent monodromy.  By $(5.4)$,
the corresponding height function $(5.5)$ is then given by the formula
$$
	\delta(1) = \frac{1}{2\pi} h(s) 1'			\tag{5.7}
$$
where $\delta$ denotes the section of $\V\otimes\V^*$ defined by the 
pointwise application of the splitting $(2.16)$ to the fibers of $\V$, 
and $1$ and $1'$ denote the generators of 
$Gr^W_0(\V)\cong\Z(0)\otimes{\Cal O}_{S'}$ and 
$Gr^W_{-2}(\V)\cong\Z(1)\otimes{\Cal O}_{S'}$ respectively.  

\par As usual, for the purpose of calculating the asymptotic behavior of
$(5.5)$ near $p$, we replace $\V$ by the corresponding period map $F:U\to\M$
obtained restricting $\V$ to a deleted neighborhood $\Delta^*$ of $p$.  Using 
the nilpotent orbit theorem discussed in \S 3, we can then replace $F(z)$ 
by the corresponding nilpotent orbit 
$$
	\theta(z) = e^{zN}.F_{\infty}				\tag{5.8}
$$
since we are only interested in calculating the leading order terms of
$(5.5)$.  Invoking Theorem $(4.2)$, we can then calculate the asymptotic
behavior of $h(s)$ modulo terms that remain bounded as $s\to 0$ (recall
$s=e^{2\pi iz}$) by replacing $\theta(z)$ by the corresponding split
orbit $\hat\theta(z) = e^{zN}.\hat F_{\infty}$.

\par Indeed, by Corollary $(4.3)$, for any admissible period map $F(z)$
of type $(\text{II})$ with unipotent monodromy, the corresponding
gradings $Y_{(F(z),W)}$ and $Y_{(\hat\theta(z),W)}$ are related by an
equation of the form
$$
	Y_{(F(z),W)} = Y_{(\hat\theta(z),W)} + \epsilon(z)	\tag{5.9}
$$
where $\epsilon(z)$ is a real analytic function which remains bounded as
$y=\text{Im}(z)\to\infty$ and $x=\text{Re}(z)$ restricted to any finite 
subinterval of $\R$.  Moreover,  by Theorem $(3.16d)$, 
$$
	Y_{(\hat\theta(z),W)} = e^{zN}.Y = Y + 2zN_{-2}		\tag{5.10}
$$
where $Y$ is a real grading of $W$, and hence 
$$
	\delta_{(\hat\theta(z),W)} = yN_{-2}			\tag{5.11}
$$
since
$$	
	Y_{(F,W)} - \overline{Y}_{(F,W)} = 4i\delta_{(F,W)}	\tag{5.12}
$$
for any mixed Hodge structure of type $(\text{II})$.  Therefore, by equation 
$(5.9)$--$(5.12)$, 
$$
	\delta_{(F(z),W)} 
	= yN_{-2} + \frac{1}{2}\Im(\epsilon(z))			\tag{5.13}
$$
Inserting equation $(5.13)$ into $(5.7)$, it then follows that, near $s=0$,
$$
	h(s) = -\mu\log|s| + \eta(s)				\tag{5.14}
$$
where $N_{-2}(1) = \mu\, 1'$ and $\eta(s)$ is a real analytic function
which remains bounded as $s\to 0$.

\remark{Remark} More generally, it follows from $(5.9)$--$(5.13)$ that if
$h_{\V}(s)$ denotes the height function $(5.1)$ attached to an admissible 
variation $\V\to\Delta^*$ of type $(\text{II})$ with unipotent monodromy,
then 
$$
	h_{\V}(s) = -\mu\log(s) + \eta(s)
$$
where $\mu=||N_{-2}||_{F_o}$ denotes the norm of $N_{-2}$ with respect to
the base point $F_o\in\M$ defined in Corollary $(4.3)$, and $\eta(s)$ is
once again a real-valued analytic function which remains bounded as $s\to 0$.
\endremark
\vskip 3pt

\par Now, according to the above recipe, in order to calculate the asymptotic
behavior of the height paring $\langle Z_s,W_s\rangle$, it would seem  that 
one must compute $N$, $W$, and $F_{\infty}$, along with the corresponding 
splittings and gradings.  This is not necessary. Indeed, these auxiliary
object appear in equation $(5.14)$ only via the decomposition
$$
	N = N_0 + N_{-2}					\tag{5.15}
$$
which can computed directly from the pair $(N,W)$ as follows:  Let $Y$ be
the grading appearing in $(5.10)$, relative to which $N$ decomposes as
$(5.15)$ according to the eigenvalues of $ad\,Y$, and $Y'$ be any other
grading of $W$.  Then, since $Lie_{-1}(W)$ acts transitively on the set of 
all gradings of $W$, 
$$
	Y' = Y + \a_{-1} + \a_{-2}				\tag{5.16}
$$
where $\a_j$ belongs to the $j$ eigenspace $E_j(\ad\, Y)$ of $\ad\, Y$.
Furthermore, because $N_0$ acts trivially on $E_0(Y)$ and $E_{-2}(Y)$,
$$
	[N_0,\a_{-2}] = 0					\tag{5.17}
$$
Therefore,
$$
	[Y',N] = [Y + \a_{-1} + \a_{-2}, N_0 + N_{-2}]
	       = -2N_{-2} + [\a_{-1},N_0]			\tag{5.18}
$$
by virtue of equation $(5.17)$ and the short length of $W$, which forces
both $[\a_{-1},N_{-2}]$ and $[\a_{-2},N_{-2}] = 0$.  Accordingly, if $Y'$ 
is any grading of $W$ such that $[Y',N]$ lowers $W$ by 2 
(i.e\. $[\a_{-1},N_0] = 0$) then $N_{-2} = -\half[Y',N]$.  Thus, in summary, 
we obtain the following result:

\proclaim{Theorem (5.19)} Let $h(s)$ denote the height function $(5.5)$ 
attached to flat family of algebraic cycles $Z_s$, $W_s\subseteq X_s$ over a 
smooth curve $S$.  Let $p$ be a point at which the corresponding 
variation $\V$ defined by equation $(5.6)$ degenerates with unipotent 
monodromy.  Let $N$ denote the local monodromy of $\V$ about $p$, and $Y'$ 
be any grading of the weight filtration $W$ of $\V$ such that $[Y',N]$
lowers $W$ by $2$.  Define $N_{-2} = -\frac{1}{2}[Y',N]$.  Then, near $s=0$,
$$
      h(s) = -\mu\log|s| + \eta(s)
$$
where $N_{-2}(1) = \mu\, 1'$ and $\eta(s)$ is a real analytic function which 
remains bounded as $s\to 0$.
\endproclaim
\demo{Proof} It remains only to justify $(5.9)$, from which Theorem $(5.19)$
then follows from equations $(5.10)$--$(5.17)$ and accompanying arguments.
To verify $(5.9)$, recall that by Corollary $(4.3)$, near the given 
puncture, the period map $F(z)$ of the variation $(5.6)$ assumes the form
$$
	F(z) = e^{xN}g(y)e^{iyN_{-2}}y^{-H/2}e^{\g(z)}.F_o	\tag{5.20}
$$
where $H\in\gg_{\R}$ commutes with the grading $Y=Y_{(F_o,W)}$ appearing
in equation $(5.10)$, and $\g(z)$ is a real analytic, $\hh$-valued 
function which is of order $y^{\b}e^{-2\pi y}$ as $y\to\infty$ and $x$  
restricted to a finite subinterval of $\R$.  Therefore,
$$
\aligned
	Y_{(F(z),W)} 
	&= e^{xN}g(y)e^{iyN_{-2}}y^{-H/2}e^{\g(z)}.Y_{(F_o,W)}		\\
	&= e^{xN}g(y)e^{iyN_{-2}}(y^{-H/2}e^{\g(z)}y^{H/2})
	   y^{-H/2}.Y_{(F_o,W)}   \\
	&= e^{xN}g(y)e^{iyN_{-2}}e^{\g_1(z)}.Y_{(F_o,W)}  
\endaligned\tag{5.21}
$$
where $\g_1(z) = \Ad(y^{-H/2})\g(z)$ is a real analytic function of order 
$y^{\b'}e^{-2\pi y}$ for some constant $\b'\in\R$.  Accordingly,
$$
	e^{\g_1(z)}.Y_{(F_o,W)} 
	= e^{\g_1(z)}.Y 
	= Y + \g_2(z)                                          \tag{5.22}
$$
where $\g_2(z)$ is again of order $y^{\b'}e^{-2\pi y}$.  Inserting $(5.22)$ 
into $(5.21)$, it then follows that 
$$
\aligned
	Y_{(F(z),W)} 
	&= e^{xN}g(y)e^{iyN_{-2}}.(Y+\g_2(z))				\\
        &= e^{xN}g(y).(Y + 2iy N_{-2} + \g_3(z))			\\
	&= (e^{xN}g(y)e^{-xN})e^{xN}.(Y + 2iy N_{-2} + \g_3(z))		\\
        &= (e^{xN}g(y)e^{-xN}).(Y + 2z N_{-2} + \g_4(z))		\\
	&= (e^{xN}g(y)e^{-xN}).(Y + \g_4(z)) 
	  +(e^{xN}g(y)e^{-xN}).(2z N_{-2}) 				
\endaligned							\tag{5.23}
$$
where, for some constant $\b''\in\R$, $\g_3(z)$ and $\g_4(z)$ are real
analytic functions of order $y^{\b''}e^{-2\pi y}$ as $y\to\infty$ with $x$ 
restricted to a finite subinterval of $\R$.  Moreover, by Theorem $(4.2)$, 
the function $g(y)$ admits a convergent series expansion near $y=\infty$ of 
the form
$$
	g(y) = e^{\zeta}(1+g_1 y^{-1} + g_2 y^{-2} + \cdots)	\tag{5.24}
$$
where $\zeta, g_1, g_2,\dots\in\ker(\ad\, N_{-2})$, and hence
$$
	g(y).N_{-2} = \Ad(g(y))N_{-2} = N_{-2}			\tag{5.25}
$$
Therefore, 
$$
	(e^{xN}g(y)e^{-xN}).(2zN_{-2}) = 2zN_{-2}		\tag{5.26}
$$
since $N=N_0 + N_{-2}$ and $[N_0,N_{-2}] = 0$.  Likewise, because of the
series expansion $(5.24)$ and the fact that  
$\zeta\in\ker(N_0)\cap\ker(N_{-2})$ by Theorem $(4.2)$, 
$$
	\lim_{y\to\infty}\, e^{xN}g(y)e^{-xN} = e^{\zeta}	\tag{5.27}
$$
independent of $x$.  Consequently,
$$
	(e^{xN}g(y)e^{-xN}).(Y + \g_4(z)) = Y + \epsilon(z)	\tag{5.28}
$$
where $\epsilon(z)$ is a real analytic function which remains bounded as
$y\to\infty$ and $x$ restricted to a finite subinterval of $\R$.  Inserting
$(5.26)$ and $(5.28)$ into $(5.23)$, it then follows that
$$
	Y_{(F(z),W)} 
	= Y + 2zN_{-2} + \epsilon(z) = Y_{(\hat\theta(z),W)} + \epsilon(z)
$$
as required.
\enddemo

\par Returning now to the general setting $(5.5)$, let $D$ be a normal
crossing divisor about which $\V$ degenerates with unipotent monodromy.  Let 
$(s_1,\dots,s_m)$ be local coordinates on $\overline{S'}$ relative to which 
$D$ assumes the form 
$
	s_1\cdots s_r = 0
$
and $f:\Delta\to \overline{S'}$ be a holomorphic map of the form 
$$
	f(t) = (t^{a_1}f_1(t),\dots,t^{a_m}f_m(t))		\tag{5.29}
$$
where $a_1,\dots, a_m$ are nonnegative integers and $f_1,\dots,f_m$ are
nonvanishing holomorphic functions on $\Delta$. 

\par Let $N_j$ denote the monodromy logarithm of $\V$ about $s_j=0$ and
$N$ denote the monodromy of $f^*(\V)$ about $t=0$.  Then,
$
	N = \sum_{j=1}^r\, a_j N_j
$
and hence
$$
	f^*(h)(t) = -\mu\log|t| + \eta(t)			\tag{5.30}
$$
where $\mu(1') = -\half[Y',N](1)$ for any grading $Y'$ of $W$ such that
$[Y',N](W_0)\subseteq W_{-2}$ and $\eta(t)$ is a real analytic function which 
remains bounded as $t\to 0$.  In particular, as consequence of the above 
remarks, the function 
$$
	\mu = \mu_{a_1,\dots,a_r}
$$ 
defined by equation $(5.30)$ is a homogeneous function of degree 1 in 
$a_1,\dots, a_r$.  As such, the simplest possible asymptotic behavior
that $h$ can exhibit as $s$ approaches $D$ along various curves of the 
form $(5.29)$ is for $\mu$ to be a linear function of $a_1,\dots,a_r$.
In this case, we shall say that $h(s)$ has no jumps along $D$.

\par By Theorem $(5.19)$, a sufficient condition for $h(s)$ to have no jumps 
along $D$ is the existence of a grading $Y$ of $W$ such that
$[Y,N_j](W_0)\subseteq W_{-2}$ for all $j$.  Indeed, in this case 
$$
   \text{$\mu_{a_1,\dots,a_r}(1') = -\half[Y,\sum_j\, a_j N_j](1)$}
								\tag{5.31}
$$
The next result gives a sufficient condition for the existence of such a
grading $Y$ which depends only on the monodromy of the local system 
$$
	Gr^{\W}_{-1}(\V_{\Z}) = [R^{2d+1}_{\pi*}(\Z(d))]^*      \tag{5.32}
$$
defined by the morphism $\pi:X\to S$, and not the particular choice of 
flat cycles $Z$ and $W\subseteq X$.

\proclaim{Theorem (5.33)} Let $W'$ denote the (unshifted) monodromy weight
filtration of $Gr^{\W}_{-1}(\V)$ about $D$.  Suppose that
$$
	\dim Gr^{W'}_{-1} = \dim Gr^{W'}_{-3}
$$
Then, the exists a grading $Y$ of $W$ such that 
$[Y,N_j](W_0)\subseteq W_{-2}$ for all $j$, and hence the corresponding
height function $h(s)$ has jumps along $D$.
\endproclaim
\demo{Proof} We begin by recalling the definition of $W'$:  Let $N'_j$ denote
the monodromy logarithm of $Gr^{\W}_{-1}(\V_{\Z})$ about $s_j=0$, and 
$$
	\Cal C' 
	= \{\, \l_1 N'_1 + \dots + \l_r N'_r \mid \l_1,\dots,\l_r>0 \,\}
$$
denote the monodromy cone of $Gr^{\W}(\V)$ about $D$.  Then, by [CK], each 
element $N'\in\Cal C'$ determines the same monodromy weight filtration 
$W'=W(N')$.

\par Let $N$ be an element of the monodromy cone
$$
	\Cal C = \{\, \l_1 N_1 + \dots +\l_r N_r \mid \l_1,\dots,\l_r>0 \,\}
$$
$\rel W = \rel W(N,W)$ and $Y$ be the grading of $W$ obtained by application 
of $(3.19)$ to $N$ and $\rel Y = Y_{(F_{\infty},\rel W)}$.  Then, relative 
to $\ad\, Y$, $N=N_0 + N_{-2}$.  Likewise, due to the short length of $W$,
each $N_j$ decomposes as
$$
	N_j = (N_j)_0 + (N_j)_{-1} + (N_j)_{-2}
$$
relative to $\ad\, Y$.  Accordingly, the condition 
$[Y,N_j](W_0)\subseteq W_{-2}$ is equivalent to the assertion that 
$(N_j)_{-1} = 0$ for each $j$. 

\par To complete the proof, let $V(k)$ denote the sum of all irreducible
submodules of $V$ of highest weight $k$ with respect to the representation
of $sl_2(\C)$ defined by $N_0$ and $H = \rel Y - Y$ [cf\. Theorem $(3.16)$].
Then, by the above remarks, it is sufficient to show that
\roster
\item"(a)" $V(1) = 0 \implies (N_j)_{-1}=0$;
\item"(b)" $\half\dim V(1) = \dim Gr^{W'}_{-1} - \dim Gr^{W'}_{-3}$.
\endroster

\par To verify $(a)$, observe that
$$
	(N_j)_{-1}\in\ker(N_0)\cap E_{-1}(\ad\, H)
$$
since $[N,N_j]=0$, $H=\rel Y - Y$ and $[\rel Y,N_j] = -2N_j$.  Therefore,
if $e_0$ is a generator of $E_0(Y)$ then
$$
	u = (N_j)_{-1}(e_0)\in\ker(N_0)				\tag{5.34}
$$
because $\rho$ acts trivially on $E_0(Y)$, and hence
$$
	N_0(u) = N_0(N_j)_{-1}(e_0) = [N_0,(N_j)_{-1}]e_0 = 0
$$
Likewise,  
$$
	u\in E_{-1}(H)						\tag{5.35}
$$
since
$$
	H(u) = H(N_j)_{-1}(e_0) = [H,(N_j)_{-1}]e_0 = -(N_j)_{-1}(e_0) = -u
$$
Combining $(5.34)$ and $(5.35)$, it then follows that $u\in V(1)$.

\par Similarly, since $\rho$ acts trivially on $E_{-2}(Y)$, if 
$v\in E_{\ell}(H)\cap E_{-1}(Y)$ and $(N_j)_{-1}(v)$ is non-zero
then $\ell = 1$ since
$$
	-(N_j)_{-1}(v) = [H,(N_j)_{-1}]v = -(N_j)_{-1}H(v)
	= -\ell (N_j)_{-1}(v)
$$
Furthermore, $(N_j)_{-1}(v)\neq 0$ implies that $v\not\in\text{Im}(N_0)$
since 
$$
	v = N_0(v') \implies 
	(N_j)_{-1}(v) = (N_j)_{-1} N_0(v') = [N_0,(N_j)_{-1}]v' = 0
$$
and hence $v\in V(1)$.  Thus, $V(1)=0 \implies (N_j)_{-1} = 0$.

\par To verify $(b)$, let $k>0$ be an odd integer.  Then, 
$E_{-k}(H)\cong Gr^{W'}_{-k}$.  Indeed, by definition, $\rel W$ induces 
the shifted monodromy weight filtration $W'[1]$ on $Gr^W_{-1}$.  Accordingly, 
$H=\rel Y - Y$ induces a grading of $W'$ on $Gr^W_{-1}$, and hence
$$
	E_{-k}(H) = E_{-k}(H)\cap E_{-1}(Y) \cong Gr^{W'}_{-k}
$$
since $V$ is the direct sum of $E_{-1}(Y)\cong Gr^W_{-1}$ and the 
trivial submodules $E_0(Y)$ and $E_{-2}(Y)$.  On the other hand,
$$
	\dim E_{-k}(H) 
	= \sum_{w\geq k,\, w\equiv k \mod 2}\, \frac{1}{w+1}\dim V(w)
$$
Therefore, the equality $\dim Gr^{W'}_{-1} = \dim Gr^{W'}_{-3}$ implies that
$$
	\half\dim V(1) = \dim E_{-1}(H) - \dim E_{-3}(H)
		       = \dim Gr^{W'}_{-1} - \dim Gr^{W'}_{-3} = 0
$$
\enddemo

\proclaim{Corollary 5.36} If the local monodromy of $Gr^{\W}_{-1}(\V_{\Z})$
about $D$ is trivial then the corresponding height function $(5.5)$ has 
no jumps along $D$.
\endproclaim

\remark{Remark} A special case of $(5.36)$ is when the cycles $Z_s$ and 
$W_s$ move about in a fixed variety $X$.  This case was considered by 
Richard Hain in [H].
\endremark

\proclaim{Corollary 5.37} Let $N'\in\Cal C'$ and suppose that the induced
map
$$
	N':Gr^{W'}_{-1}\to Gr^{W'}_{-3}
$$
is injective.  Then, the associated height function $h(s)$ has no jumps
along $D$.
\endproclaim
\demo{Proof} By standard $sl_2$ theory, 
$\dim V(1) = 2\dim\ker(N':Gr^{W'}_{-1}\to Gr^{W'}_{-3})$.
\enddemo


\par To close this section, we now present two related examples which show 
that when $\dim Gr^{W'}_{-1} \neq \dim Gr^{W'}_{-3}$ the height may or may
not jump:

\example{Example 5.38} Let $V_{\Z}$ be an integral lattice of rank 4, with 
basis $\{e_0,e,f,e_{-2}\}$, and $N_1$, $N_2$ denote the endomorphisms of 
$V_{\Z}$ defined by the matrices
$$
	N_1 = \pmatrix 0 & 0 & 0 & 0 \\
		       0 & 0 & 0 & 0 \\
		       1 & 1 & 0 & 0 \\ 
		       1 & 1 & 0 & 0 \endpmatrix,\qquad
	N_2 = \pmatrix  0 &  0 & 0 & 0 \\
		        0 &  0 & 0 & 0 \\
		       -1 &  1 & 0 & 0 \\ 
		        1 & -1 & 0 & 0 \endpmatrix
$$
Then, the variation $\V\to\Delta^{*2}$ defined by the weight filtration
\roster
\item"" $W_0(V_{\Z}) = V_{\Z}$;
\item"" $W_{-1}(V_{\Z}) = \Span_{\Z}(e,f,e_{-2})$;
\item"" $W_{-2}(V_{\Z}) = \Span_{\Z}(e_{-2})$;
\item"" $W_{-3}(V_{\Z}) = 0$;
\endroster
and the period map
$$
	\varphi(s_1,s_2)
  = e^{\frac{1}{2\pi i}(\log(s_1)N_1 + \log(s_2)N_2)}.F_{\infty},
$$
defined by $N_1$, $N_2$ and the filtration
\roster
\item"" $F^{-1}_{\infty} = \Span_{\C}(e_0,e,f,e_{-2})$;
\item"" $F^0_{\infty} = \Span_{\C}(e_0,e)$;
\item"" $F^1_{\infty} = 0$;
\endroster
is admissible, and graded-polarizable.  Direct calculation shows that
the associated height function $(5.5)$ is given by the formula
$$
	h(s_1,s_2) 
	= \frac{(\log|s_1/s_2|)^2 -(\log|s_1 s_2|)^2}
		{\log|s_1 s_2|}					\tag{5.39}
$$
Setting $(s_1,s_2) = (t^{a_1},t^{a_2})$, it then follows that
$$
	\mu = \frac{4a_1 a_2}{a_1 + a_2}
$$
is a non-linear function of $(a_1,a_2)$ [i.e\. $h(s_1,s_2)$ jumps].
\endexample

\example{Example 5.40} In Example $(5.38)$, redefine
$$
	N_1 = \pmatrix 0 & 0 & 0 & 0 \\
		       0 & 0 & 0 & 0 \\
		       0 & 1 & 0 & 0 \\ 
		       1 & 1 & 0 & 0 \endpmatrix,\qquad
	N_2 = \pmatrix  0 &  0 & 0 & 0 \\
		        0 &  0 & 0 & 0 \\
		        0 &  1 & 0 & 0 \\ 
		        1 & -1 & 0 & 0 \endpmatrix
$$
Then, 
$$
	h(s_1,s_2) = -\log|s_1 s_2|				
$$
and hence $\mu = a_1 + a_2$.  Accordingly, $h(s_1,s_2)$ has no jumps along 
$D$. 
\endexample

\head \S 6.\qquad Nahm's Equation \endhead

\par Let $K$ be a compact real Lie group.  Then, Nahm's equation for $K$
is the system of ordinary differential equations given by the 
gradient flow of the 3-form
$$
	\phi(T_1,T_2,T_3) = \langle T_1,[T_2,T_3]\rangle	\tag{6.1}
$$
on $\kappa = Lie(K)$ defined by a choice of bi-invariant metric
$\langle\cdot,\cdot\rangle$ on $K$.  Equivalently, a triple of 
$\kappa$-valued functions $(T_1,T_2,T_3)$ is a solution of Nahm's
equation if and only if 
$$
	\frac{dT_i}{dy} + [T_j,T_k] = 0				\tag{6.2}
$$
for every cyclic permutation $(i\hph{,} j\hph{,} k)$ of 
$(1\hph{,} 2\hph{,} 3)$.  

\par More generally, given a complex Lie algebra $\aa$, a triple of 
$\aa$-valued functions $(T_1,T_2,T_3)$ is said to be a solution of Nahm's 
equation provided they satisfy the system of differential equations $(6.2)$. 
Solutions to Nahm's equation are related to representations of $sl_2(\C)$
as follows:  Let $\{\tau_1,\tau_2,\tau_3\}$ be a basis of 
$sl_2(\C)=su_2\otimes\C$ such
that
$$
	\tau_i = [\tau_j,\tau_k]				\tag{6.3}
$$
for every cyclic permutation $(i\hph{,} j\hph{,} k)$ of 
$(1\hph{,} 2\hph{,} 3)$ and $\rho:sl_2(\C)\to\aa$ be a Lie algebra 
homomorphism.  Then, the triple
$$
	T_i(y) = \rho(\tau_i)y^{-1}
$$
is a solution of $(6.2)$.  Conversely, given a solution $(T_1,T_2,T_3)$
of Nahm's equation which has a simple pole at $y=0$, the linear map
$\rho:sl_2(\C)\to\aa$ defined by setting 
$$
	\rho(\tau_i) = \text{Res}(T_i)
$$
is a Lie algebra homomorphism.
\vskip 3pt

\par In [S], Schmid showed that a nilpotent orbit of pure, polarized Hodge
structure gives rise to a solution 
$$
	\Phi:(a,\infty)\to \Hom(sl_2(\C),\gg_{\C}),\qquad
	\Phi(y)\tau_i = T_i(y)					\tag{6.4}
$$
of Nahm's equation.  In this section, we show that a nilpotent orbit 
$$
	\theta(z) = e^{zN}.F_{\infty}				\tag{6.5}
$$
of graded-polarized mixed Hodge structure gives rise to a solution of a
generalization of Nahm's equation which encodes how the extension data 
of $\theta(z)$ interacts with the nilpotent orbits of pure Hodge structure 
induced by $\theta(z)$ on $Gr^W$.
\vskip 3pt

\par To this end, let $\M$ be a classifying space of graded-polarized
mixed Hodge structure.  Define $\Cal D$ to be the direct sum of 
classifying spaces of pure, polarized Hodge structure onto which
$\M$ projects via the map
$$
	F\mapsto FGr^W
$$
Let $\tlY(Y)$ be the affine space consisting of all gradings $Y$ of $W$ such 
that [cf\. Theorem $(2.17)$]:
$$
	Y - \bar Y\in Lie_{-2}(W)				\tag{6.6} 
$$
and $\iota_Y$ denote the isomorphism $Gr^W\cong V$ associated to 
$Y\in\tlY(W)$.  Then:

\proclaim{Theorem 6.7} The space $\X=\Cal D\times\tlY(W)$ is a complex
manifold upon which the Lie group $\tlG$ \cite{cf\. Theorem $(2.19)$} acts
transitively by automorphisms.  Furthermore, the correspondence
$$
	F = \pi(\{H^{r,s}\},Y) \iff 
	F^p = \bigoplus_{a\geq p}\, \iota_Y(H^{a,b})		\tag{6.8}
$$
defines a $\tlG$--equivariant projection map $\pi:\tlX\to\M$ with real
analytic section  
$$
	\sigma(F) = (FGr^W,Y_{(F,W)})				\tag{6.9}
$$
\endproclaim
\demo{Proof} The only subtle point is the assertion that $\sigma$ is a
real--analytic section.  To prove this, observe that by part
$(c)$ of Theorem $(2.4)$, the grading $Y_{(F,W)}$ defined by the
$I^{p,q}$'s of $(F,W)$ takes values in $\tlY(W)$.  Consequently,
equation $(6.9)$ defines a section of $\X$.  To prove that $\sigma$
is real-analytic, recall [CKS] that
$$
	I^{p,q} = F^p\cap W_{p+q}
		  \cap( \bar F^q\cap W_{p+q} 
		  + \sum_{j>0}\, \bar F^{q-j}\cap W_{p+q-1-j})
$$
and hence the decomposition $(2.5)$ is real-analytic with respect to
the point $F\in\M$.
\enddemo

\par Next, following [S], we note that each choice of base point 
$F_o$ defines a principal bundle $P$ over $\X$ with connection $\nabla$:

\proclaim{Theorem 6.10} Let $F_o\in\M_{\R}$ and $x_o = \sigma(F_o)$.  Then,
the vector space \cite{cf\. Theorem $(2.12)$}
$$
	\hh' = (\np\oplus\lam\oplus\nn)\cap\hh
$$
is an $\Ad(\tlG^{x_o})$--invariant complement to $\hh^{x_o}$ in $\hh$,
and hence defines a connection $\nabla$ on the principal bundle
$$
	\tlG^{x_o}\to\tlG\to\tlG/\tlG^{x_o}
$$
over $\X\cong\tlG/\tlG^{x_o}$.
\endproclaim
\demo{Proof} Direct calculation shows that since $F_o\in\M_{\R}$, 
$\hh^{x_o} = \nz\cap\hh$ and hence $\hh'$ is a vector space 
complement to $\hh^{x_0}$ in $\hh$.  To see that $\hh'$ is invariant under 
the action of $\Ad(H^{x_0})$, let $h\in H^{x_o}$.  Then, $h$ preserves $F_o$ 
since
$$
	h.F_o = h.\pi(x_o) = \pi(h.x_o) = \pi(x_o) = F_o
$$
Likewise, $h = \bar h$ since $h$ acts by real automorphisms on $Gr^W$
and preserves the real grading $Y_{(F_o,W)}$.  Consequently, $h$ is 
a morphism of $(F_o,W)$ and hence preserves each summand appearing in 
the definition of $\hh'$.
\enddemo

\par Thus, by virtue of the above remarks, each choice of base point 
$F_o\in\M_{\R}$ defines a lift of $\theta(iy)$ to a function 
$h(y):(a,\infty)\to\tlG$ such that:
\roster
\item"(a)" $h(y).x_o = \sigma(\theta(iy))$;
\item"(b)" $h$ is tangent to $\nabla$.
\endroster

\proclaim{Theorem 6.11} Let $\LL$ denote the endomorphism of $\hh$ defined
by the rule:
$$
	\left.\LL\right|_{\np} = +i,\qquad
	\left.\LL\right|_{\nz} =  0,\qquad	
	\left.\LL\right|_{\nn\oplus\lam} = -i
$$
Then, the function $h(y)$ defined above satisfies the differential equation
$$
	h^{-1}(y)\frac{d}{dy}h(y) = -\LL\Ad(h^{-1}(y))N		\tag{6.12}
$$
\endproclaim
\demo{Proof} Schmid's original derivation [S, Lemma $(9.8)$] of Nahm's 
equation for nilpotent orbits of pure, polarized Hodge structure shows that
equation $(6.12)$ holds modulo $Lie_{-1}(W)$.  Consequently, it is sufficient 
to verify that equation $(6.12)$ holds modulo the subalgebra 
$\gg_{\C}^Y = Lie(\Gc^Y)$, $Y=Y_{(F_o,W)}$ since
$$
	\gg_{\C} = \gg_{\C}^Y \oplus Lie_{-1}(W)
$$
\par To this end, note that by definition 
$Y_{e^{iyN}.F_{\infty}} = \Ad(h(y))Y$.  Upon differentiating both sides of 
this equation with respect to $y$ and simplifying the result, it then follows 
that:
$$
	\Ad(h^{-1}(y))\frac{d}{dy}\, Y_{(e^{iyN}.F_{\infty},W)}
	= \left[h^{-1}(y)\frac{d}{dy}h(y),Y\right]		\tag{6.13}
$$
Therefore, if $z=x+iy$:
$$
        \Ad(h^{-1}(y))
	\frac{d}{dy} Y_{(e^{iyN}.F_{\infty},W)}
	= i\,\left.\Ad(h^{-1}(y))
             \left(\frac{\pd}{\pd z} - \frac{\pd}{\pd\bar z}\right)
              Y_{(e^{zN}.F_{\infty},W)}\right|_{z=iy}
	\tag{6.14}
$$

\par To compute $\frac{\pd}{\pd z}\,Y_{(e^{zN}.F_{\infty},W)}$ and 
$\frac{\pd}{\pd\bar z}\,Y_{(e^{zN}.F_{\infty},W)}$, we observe that
as a consequence of equation $(5.19)$ in [P2]:
$$
\aligned
   \left.\frac{\pd}{\pd w} Y_{(e^{w\xi}.F,W)}\right|_{w=0} 
   &= [\piq(\xi),Y_{(F,W)}],						\\
   \left.\frac{\pd}{\pd\bar w} Y_{(e^{w\xi}.F,W)}\right|_{w=0}
   &= [\pip(\overline{\piq(\xi)}),Y_{(F,W)}]				
\endaligned\tag{6.15}
$$
for any point $F\in\M$ and any element $\xi\in Lie(\Gc)$, where $\pip$ and 
$\piq$ denote the projection operators\footnote{In [P2], we used the 
alternative notation ${\frak t}_F = q_F$ and $\piq = \pi_q$.} with respect to 
$F$ defined in Theorem $(2.12)$.  In particular, upon setting
$F=e^{iy N}.F_{\infty}$ it then follows from equations $(6.14)$ 
and $(6.15)$ that:
$$
   \Ad(h^{-1}(y))
	  \frac{d}{dy}\, Y_{e^{iyN}.F_{\infty}}
	= i\,\Ad(h^{-1}(y))
             [\piq(N)-\pip(\overline{\piq(N)}),Y_{e^{iy N}.F_{\infty}}]
	\tag{6.16}
$$

\par On the other hand, if $\piz$ denotes projection onto $\nz$ with
respect to $F=e^{iy N}.F_{\infty}$ then
$$
	N = \pip(N) + \piz(N) + \piq(N)
$$
Consequently, since $N$ is defined over $\R$:
$$
	N = \bar N = \overline{\pip(N)} + \overline{\piz(N)} 
		     + \overline{\piq(N)}
$$
and hence $\pip(N) = \pip(\overline{\piq(N)})$.  Accordingly, 
equation $(6.16)$ may be rewritten as
$$
\aligned
   \Ad(h^{-1}(y)) &
	\frac{d}{dy}\, Y_{e^{iyN}.F_{\infty}}
	= i\,\Ad(h^{-1}(y))
             [\piq(N)-\pip(N),Y_{e^{iy N}.F_{\infty}}]			\\
	&= i[\Ad(h^{-1}(y))\{\piq(N)-\pip(N)\},
             \Ad(h^{-1}(y))Y_{e^{iy N}.F_{\infty}}]			\\
        &= i[\Ad(h^{-1}(y))\{\piq(N)-\pip(N)\},Y]
\endaligned\tag{6.17}
$$
since $Y_{(e^{iy N}.F_{\infty},W)} = Ad(h(y))Y$.  

\par By construction:
$$
	h(y).I^{p,q}_{(F_0,W)} = I^{p,q}_{(e^{iy N}.F_{\infty},W)}
								\tag{6.18}
$$
and hence
$
	\Ad(h(y)):\gg^{r,s}_{(F_0,W)}\to 
		    \gg^{r,s}_{(e^{iy N}.F_{\infty},W)}
$.
Consequently, 
$$
\aligned
	i\,\Ad(h^{-1}(y))\{\piq(N)-\pip(N)\} 
        &= i\,\hat\piq(\Ad(h^{-1}(y))N)-i\,\hat\pip(\Ad(h^{-1}(y)N))     \\
	&= -\LL\,\Ad(h^{-1}(y))N \mod Lie(\Gc^Y) 
\endaligned
$$
where $\hat\piq$ and $\hat\pip$ denote projection with respect to 
$F_o\in\M_{\R}$.  Therefore, by equation $(6.17)$,
$$
   \Ad(h^{-1}(y))
	\frac{d}{dy}\, Y_{(e^{iyN}.F_{\infty},W)}
        = [-\LL\,\Ad(h^{-1}(y))N,Y]
   								\tag{6.19}
$$
Accordingly, upon comparing equation $(6.19)$ with equation $(6.13)$, it 
then follows that
$$
	[-\LL\,\Ad(h^{-1}(y))N,Y] = [h^{-1}(y)\frac{d}{dy}h(y),Y]
$$	
and hence
$-{\LL}\, \Ad(h^{-1}(y))N = h^{-1}(y)\frac{d}{dy}h(y)$ mod $\gg_{\C}^Y$
as required.
\enddemo

\example{Example 6.20} Let $\theta(z) = e^{zN}.\hat F$ be a split orbit.  
Then, the function 
$$
	h(y) = e^{iyN}e^{-iyN_0}y^{-H/2}			
$$
[cf\. Theorem $(3.16)$ for notation] is a solution of equation $(6.12)$ with 
respect to the base point $F_o = e^{iN_0}.\hat F\in\M_{\R}$.  

\par To prove this, equip $sl_2(\C)$ with the standard Hodge structure
$(3.11)$ and $\gg_{\C}$ with the usual mixed Hodge structure induced by
$(F_o,W)$.  Then, as a consequence of Theorem $(3.13)$ and the fact
[Theorem $(3.16)$, part $(c)$] that $e^{zN_0}.\hat F$ is and 
$\text{\rm SL}_2$-orbit with data $(F_o,\psi_*=\rho)$, the representation
$$
	\rho:sl_2(\C)\to\gg_{\C}				\tag{6.21}
$$
defined in Theorem $(3.16)$ is morphism of Hodge structure.

\par By direct calculation:
$$
\aligned
	h^{-1}\frac{dh}{dy} 
	&= -\frac{H}{2y} 
	   + i\Ad(y^{H/2})Ad(e^{iyN_0})(\sum_{k\geq 2}\, N_{-k})	\\
	\Ad(h^{-1}(y)) N &= \frac{N_0}{y}
	   + \Ad(y^{H/2})Ad(e^{iyN_0})(\sum_{k\geq 2}\, N_{-k})
\endaligned							\tag{6.22}
$$
Similarly, a small computation in $sl_2(\C)$ shows that the basis
$(1.6)$ satisfies the Hodge conditions:
$$
	\xx^-\in sl_2(\C)^{-1,1},\quad		
	\zz\in sl_2(\C)^{0,0},\quad
	\xx^+\in sl_2(\C)^{1,-1}				\tag{6.23}
$$
Therefore, since $\rho$ is a morphism of Hodge structures, the image
$(X^+,Z,X^-)$ of the basis $(1.6)$ under $\rho$ satisfy the analogous
conditions 
$$
	X^- \in\gg^{-1,1},\quad
	Z   \in\gg^{0,0},\quad	
	X^+ \in\gg^{1,-1}					\tag{6.24}
$$
at $(F_o,W)$.  Comparing $(1.6)$ and $(3.15)$, it then follows that
$$
\gathered
	N_0 = \frac{1}{2i}(X^+ - X^- + Z),\qquad
	N_0^+ = \frac{1}{2i}(X^+ - X^- - Z)			\\
	H = (X^+ + X^-)			
\endgathered							\tag{6.25}
$$
Consequently,
$$
	\LL(N_0) = \frac{1}{2i}\LL(X^+ - X^- + Z)
		 = \frac{1}{2i}(iX^+ + iX^-) = \half H		\tag{6.26}
$$

\par To continue, we now recall that by [D4] [KP]
$$
	(\ad N_0)^j N_{-k}\in\lam_{(F_o,W)}			\tag{6.27}
$$
and hence the function 
$
	\Ad(y^{H/2})Ad(e^{iyN_0})(\sum_{k\geq 2}\, N_{-k})
$
takes values in $\lam_{(F_o,W)}$.  Therefore, by equations $(6.22)$
and $(6.26)$:
$$
\aligned 
	-\LL\Ad(H^{-1}(y)) N 
	&= -\frac{\LL(N_0)}{y}
	   - \LL\Ad(y^{H/2})Ad(e^{iyN_0})(\sum_{k\geq 2}\, N_{-k})	\\
	&= -\frac{H}{2y}
	   +i\Ad(y^{H/2})Ad(e^{iyN_0})(\sum_{k\geq 2}\, N_{-k})         
	 = h^{-1}\frac{dh}{dy} 
\endaligned
$$
\endexample

%

\par To relate equation $(6.12)$ with Nahm's equation, we now decompose
$$
	\b(y) = \Ad(h^{-1}(y))N					\tag{6.28}
$$
according to its Hodge components with respect to $(F_o,W)$.  To this
end, observe that as a consequence of equation $(6.18)$, the Hodge
decomposition of $\b(y)$ with respect to $(F_o,W)$ has the same form
as the Hodge decomposition of $N$ with respect to $(e^{zN}.F_{\infty},W)$.
Therefore, by the next lemma, the Hodge decomposition of $\b(y)$ with
respect to $(F_o,W)$ is of the form
$$
	\b(y) =   \b^{1,-1}(y) + \b^{0,0}(y) + \b^{-1,1}(y) 
		+ \b_+(y) + \b_-(y)				\tag{6.29}
$$
where 
$$
	\b_+(y)   = \sum_{k>0}\, \b^{0,-k}(y),\qquad	
	\b_-(y) = \sum_{k>0}\, \b^{-1,1-k}(y)			\tag{6.30}
$$

\proclaim{Lemma 6.31} Let $e^{zN}.F$ be a nilpotent orbit.  Then,
with respect to $(e^{zN}.F,W)$, the Hodge decomposition of $N$ 
assumes the form:
$$  
	N = N^{-1,1} + N^{0,0} + N^{1,-1} 
	    + \left(\sum_{k>0} N^{-1,1-k}\right) 
	    + \left(\sum_{k>0} N^{0,-k}\right)			\tag{6.32}
$$
\endproclaim
\demo{Proof} The fact the $N$ is horizontal at $e^{zN}.F$ implies that
$$
	N = N^{-1,1} + \sum_{k>0}\, N^{-1,1-k} \mod 
	    \bigoplus_{r\geq 0}\,\gg^{r,s}			\tag{6.33}
$$
Define
$$
	N_{-k} = \bigoplus_{r+s = -k}\, N^{r,s}
$$
Then, the horizontality $(6.33)$ of $N$ coupled with the fact that 
$N = \bar N$ implies that
$$
	N_0 = N^{-1,1} + N^{0,0} + N^{1,-1}			\tag{6.34}
$$

\par Suppose that $(6.32)$ is false and let $k$ be the smallest integer such
that $N_{-k}$ violates $(6.32)$.  By $(6.34)$, $k>0$.  As such, by
equation $(6.33)$
$$
	N_{-k} = N^{-1,1-k} + N^{0,-k} + N^{p,-p-k} + \cdots	
$$
for some integer $p>0$.  By, Theorem $(2.4)$:
$$
	\overline{\gg^{r,s}} 
	= \gg^{s,r} \mod \bigoplus_{a<s,b<r}\,\gg^{a,b}		\tag{6.35}
$$
Accordingly, $\overline{N^{p,-k-p}}$ is of Hodge type $(-k-p,p)$ modulo
lower order terms.  Consequently, since $N=\bar N$ and elements of type 
$(-k-p,p)$ are not horizontal, $\overline{N^{p,-k-p}}$ must be annihilated 
by part of the fallout of the complex conjugate of some Hodge component 
$N^{r,s}$ with $r+s>-k$.  On the other hand, by the definition of $k$, all 
such components $N^{r,s}$ satisfy $(6.32)$.  Therefore, by equation $(6.35)$, 
there is no way for $\overline{N^{r,s}}$ to annihilate $\overline{N^{p,-k-p}}$
since $p>0$.
\enddemo

\par Following [S] [CKS], define
$$
	\a(y) = -2h^{-1}(y)\frac{dh}{dy}		\tag{6.36}
$$
Then, by virtue of equation $(6.12)$,
$$
	\a(y) = \a^{1,-1}(y) + \a^{-1,1}(y) + \a_+(y) + \a_-(y)
$$
where 
$$
\aligned
	\a^{1,-1} &= 2i\b^{1,-1},	\\
	\a^{-1,1} &= -2i\b^{-1,1} 
\endaligned\qquad
\aligned
	\a_+ &= 2i\b_+			\\
	\a_- &= -2i\b_-
\endaligned							\tag{6.37}
$$
On the other hand, differentiation of equation $(6.28)$ shows that 
$$
	-2\frac{d\b}{dy} = [\b(y),\a(y)]  			\tag{6.38}
$$
Inserting equation $(6.37)$ into $(6.38)$ and taking Hodge components,
we then obtain the following result:

\proclaim{Theorem 6.39} Let $h(y)$ be a solution to equation $(6.12)$.
Then, 
$$
         \frac{d}{dy}\b_0(y) = -[\b_0(y),\LL \b_0(y)],\qquad
         \b_0(y) = \sum_{r+s=0}\, \b^{r,s}(y)                      \tag{6.40}
$$
and 
$$
           \frac{d}{dy}\pmatrix \b_- \\ \b_+ \endpmatrix
           = i\pmatrix   \ad\, \b^{0,0}  & -2\,\ad\, \b^{-1,1} \\
                     2\,\ad\, \b^{1,-1} & -\ad\, \b^{0,0}  \endpmatrix
             \pmatrix \b_- \\ \b_+ \endpmatrix
	     +2i\pmatrix [\b_+,\b_-] \\ 0 \endpmatrix 		\tag{6.41}
$$
\endproclaim

\par In particular, as a consequence of equation $(6.40)$, we obtain
the following relationship between nilpotent orbits and solutions
to Nahm's equation:

\proclaim{Corollary 6.42} Let $h(y)$ be a solution of equation $(6.12)$,
and 
$$
	X^-(y) = -2i\b^{-1,1}(y),\quad
	Z(y) = 2i\b^{0,0}(y),\quad				 
	X^+(y) = 2i\b^{1,-1}(y)					\tag{6.43}
$$
The function $\Phi:(a,\infty)\to\text{\rm Hom}(sl_2(\C),\gg_{\C})$ 
defined by setting 
$$
	\Phi(y) \xx^+ = X^+(y),\quad
	\Phi(y) \zz = Z(y),\quad				
	\Phi(y) \xx^- = X^-(y)				\tag{6.44}
$$
is a solution \cite{$(6.4)$} of Nahm's equation.
\endproclaim
\demo{Proof} The assertion that $\Phi$ is a solution to Nahm's equation
is equivalent to the system of equations:
$$
\gathered
	-2\frac{dX^+}{dy} = [Z(y),X^+(y)],\qquad
	2\frac{dX^-}{dy}  = [Z(y),X^-(y)]			\\	
	-\frac{dZ}{dy}    = [X^+(y),X^-(y)]		
\endgathered							\tag{6.45}
$$
To verify that the triple $(6.43)$ satisfies equation $(6.45)$, 
one simply expands out equation $(6.40)$ in terms of the Hodge components
of $\b_0$.
\enddemo

\par The remaining Hodge components of $\a$ and $\b$ determine the 
extension data $\theta(z)$.  To relate these components to solution
of Nahm's equation $(6.44)$, let
$$
	A = \pmatrix \frac{1}{2}\ad\,Z(y) & -\ad\, X^-(y) \\
                    -\ad\, X^+(y) & -\frac{1}{2}\ad\, Z(y) \endpmatrix
		    						\tag{6.46}
$$
and define 
$$
	\tau_{-k} = \sum_{r>0,s>0,r+s=k}\, [\a^{0,r},\a^{-1,1-s}]
								\tag{6.47}
$$
Then, equation $(6.41)$ is equivalent to the hierarchy of differential
equations:
$$
	\frac{d}{dy}\pmatrix \a^{-1,1-k} \\ \a^{0,-k} \endpmatrix 
	= A\pmatrix \a^{-1,1-k} \\ \a^{0,-k} \endpmatrix 
	  + \pmatrix \tau_{-k} \\ 0 \endpmatrix,\qquad
	  k=1,2,... 						\tag{6.48}
$$
Accordingly, equation $(6.48)$ can be viewed as a system of equations
relating the evolution of the extension data of $\theta(z)$ to the 
nilpotent orbits of pure Hodge structure induced by $\theta(z)$ on 
$Gr^W$.

\head \S 7.\quad Nilpotent Orbits of Pure Hodge Structure \endhead

\par The relation between nilpotent orbits and solutions
of the generalized Nahm's equation presented in Theorem $(6.11)$
can be inverted as follows:

\proclaim{Theorem 7.1} Let $F_o\in\M_{\R}$, and suppose that $\b(y)$ is an
$\hh$-valued function which satisfies the Lax equation
$$
	\frac{d\b}{dy} = -[\b(y),\LL\b(y)]		\tag{7.2}
$$
Then, there exists an $\hh$-valued function $h(y)$, an element 
$\tilde N\in\hh$ and a point $\tilde F\in\check\M$ such that 
\roster
\item"(a)" $h^{-1}(y)\frac{dh}{dy} = -\LL\b(y)$, 
	   $\b(y) = \Ad(h^{-1}(y))\tilde N$;
\item"(b)" $h(y).F_o = e^{iy\tilde N}.\tilde F$.
\endroster
\endproclaim
\demo{Proof} The differential equation 
$$
	h^{-1}(y)\frac{dh}{dy} = -\LL\b(y)		\tag{7.3}
$$
completely determines $h(y)$ up to a choice of initial value $h_o\in\tlG$.
Likewise, by virtue of equations $(7.2)$ and $(7.3)$,
$$
	\Ad(h^{-1}(y))\frac{d}{dy}\Ad(h(y))\b(y) = 0	
$$
Therefore, $\b(y) = \Ad(h^{-1}(y))\tilde N$ for some fixed element 
$\tilde N\in\hh$.  Similarly, by virtue of equation $(7.3)$,
$$
\aligned
	h^{-1}(y)e^{iy\tilde N}\frac{d}{dy} e^{-iy\tilde N}h(y) 
	&= h^{-1}(y)e^{iy\tilde N}\left(-i\tilde N e^{-iy\tilde N}h(y) 
	            + e^{iy\tilde N}\frac{dh}{dy}\right)		\\
	&= -i\Ad(h^{-1}(y))\tilde N + h^{-1}\frac{dh}{dy}		\\
	&= -i\b(y) - \LL\b(y) \in\gg_{\C}^{F_o}
\endaligned
$$
Accordingly, 
$$
	e^{-iy\tilde N}h(y) = g_{\C}f(y)
$$ 
for some $\Gc^{F_o}$-valued function $f(y)$ and some fixed element 
$g_{\C}\in\Gc$.  Thus,
$$
	h(y).F_o 
	= e^{iy\tilde N}g_{\C}f(y).F_o
	= e^{iy\tilde N}.\tilde F
$$
where $\tilde F = g_{\C}.F_o$.
\enddemo

\remark{Remark} In order for $e^{z\tilde N}.\tilde F$ to be a proper nilpotent
orbit in the sense of Definition $(3.8)$, $\tilde N$ must be real and $\b(y)$
must be horizontal with respect to $F_o$.  In this case, we can then 
introduce a spectral parameter into equation $(7.2)$ by simply replacing 
$\b(y)$ by 
$
	\b_{\l}(y) = \sum_{p,q}\, \l^p \b^{p,q}(y)
$. 
\endremark
\vskip 3pt

\par In \S 6 of [CKS], Cattani, Kaplan and Schmid proved the 
$\text{\rm SL}_2$-orbit theorem for nilpotent orbits of pure Hodge structure 
$\theta(z) = e^{zN}.F$ by constructing a series solution
$$
	\b(y) = \sum_{n\geq 0}\, \b_n y^{-1-n/2}
$$
of equation $(7.2)$ such that $(\tilde N,\tilde F) = (N,F)$.  
In this section, we summarize this approach in some detail in preparation 
for the proof Theorem $(4.2)$ presented in \S 8--9.
\vskip 3pt

\par To this end, let $\aa$ be a complex Lie algebra and $\frak U$ be a
representation of $sl_2(\C)$.  Then, contraction against the Casimir
element
$$
	\Omega = 2\xx^+\xx^- + 2\xx^-\xx^+ + \zz^2		\tag{7.4}
$$
of $sl_2(\C)$ defines a pairing
$$
	Q:\Hom(sl_2(\C),\aa)\otimes\Hom(\UU,\aa)
	\to\Hom(\UU,\aa)				\tag{7.5}
$$
via the rule 
$$
	Q(A,B)(u) = 2[A(\xx^+),B(\xx^-.u)] + 2[A(\xx^-),B(\xx^+.u)]
		    + [A(\zz),B(\zz.u)]				\tag{7.6}
$$
Furthermore, a short calculation shows that, relative to the adjoint
representation $\frak U$ of $sl_2(\C)$, Nahm's equation $(6.2)$ is
equivalent to the differential equation 
$$
	-8\frac{d\Phi}{dy} = Q(\Phi,\Phi)			\tag{7.7}
$$

\par Following [CKS], suppose that $\Phi$ has a convergent series expansion
about infinity of the form
$$
	\Phi = \sum_{n\geq 0}\, \Phi_n y^{-1-n/2}
$$
and let $Q = 8Q_o$.  Then, equation $(7.7)$ is equivalent to the recursion
relations
$$
	\Phi_0 = Q_o(\Phi_0,\Phi_0)				\tag{7.8}
$$
and 
$$
	(1+n/2)\Phi_n - 2Q_o(\Phi_0,\Phi_n) 
	= \sum_{0<k<n}\, Q_o(\Phi_k,\Phi_{n-k}),\qquad n>0	\tag{7.9}
$$
Equation $(7.8)$ implies that $\Phi_0$ is either zero or an embedding
of $sl_2(\C)$ in $\gg_{\C}$.  If $\Phi_0 = 0$ then $\Phi_n =0$ for all
$n$ by induction.  If $\Phi_0 \neq 0$ then a short calculation 
[CKS:6.14] shows that
$$
	Q_o(\Phi_0,T) = \frac{1}{16}(\ell(\Omega) - \Omega T + 8T)
$$
where $\ell(\Omega)T$ and $\Omega T$ respectively denote the left and
diagonal action of the Casimir element $(7.4)$ on 
$T\in\Hom(sl_2(\C),\gg_{\C}) \cong \gg_{\C}\otimes sl_2(\C)^*$.

\par To continue, we now recall [CKS:6.18] that relative to the $sl_2$ 
module structure induced on $\gg_{\C}$ by $\Phi_0$, we can decompose
$$
	\Hom(sl_2(\C),\gg_{\C}) 
	= \sum_{r\geq 0}\sum_{\epsilon=-1}^1\,
	  \Hom(sl_2(\C),\gg(r))^{\epsilon}		\tag{7.10}
$$
where $\Hom(sl_2(\C),\gg(r))^{\epsilon}$ is the isotypical component of 
consisting of the span of all irreducible submodules of 
$\gg_{\C}\otimes sl_2(\C)^*$ which are of highest weight $r$ with 
respect to the left module structure and highest weight $r+2\epsilon$
with respect to the diagonal structure.

\par Relative to the bigrading $(7.10)$, the recursion relation $(7.9)$
reduces to the equation [CKS:6.20]
$$
	(n + \epsilon^2 + \epsilon(r+1))\Phi_n^{r,\epsilon}
	= 2 \sum_{0<k<n}\, Q_o(\Phi_k,\Psi_{n-k})^{r,\epsilon}	\tag{7.11}
$$
Therefore, subject to the compatibility condition
$$
	\sum_{0<k<n}\, Q_o(\Phi_k,\Psi_{n-k})^{n,-1} = 0	\tag{7.12}
$$
equation $(7.11)$ completely determines every component $\Phi_n^{r,\epsilon}$
except $\Phi_n^{n,-1}$ in terms of $\Phi_0,\dots,\Phi_{n-1}$.  The 
verification of the compatibility condition $(7.12)$ in turn reduces a
standard weight argument (cf\. [CKS:6.21]).  Thus, given a collection
of elements
$$
	T_n\in\Hom(sl_2(\C),\gg(n))^{-1}			\tag{7.13}
$$
there exists a unique series solution $\Phi$ of equation $(7.7)$ such that
\roster
\item"(a)" $\Phi_n\in\oplus_{r\leq n,\, r\equiv n \mod 2}\, 
             Hom(sl_2(\C),\gg(r))$;
\item"(b)" $\Phi_n^{n,-1} = T_n$;
\item"(c)" $\Phi_n^{n,0} = \Phi_n^{n,1} = 0$.
\endroster
In particular, $\Phi_1 = 0$ since it must highest weight $-1$ with respect
to the diagonal action of $sl_2(\C)$.  

\par Imposing the condition that $\Phi$ should be horizontal and map
$sl_2(\R)$ into $\hh = \gg_{\R}$, it then follows that each $T_n$ must also
be a morphism of Hodge structure with respect to the standard Hodge
structure on $sl_2(\C)$ defined in \S 3 and pure Hodge structure
$$
	\gg_{\C} = \bigoplus_p\, \gg^{p,-p}
$$
induced by $F_o$ on $\gg_{\C}$.  Accordingly, since  $\gg_{\C}$ is the Lie 
algebra of a linear Lie group $\Gc$, the equation 
$$
	h^{-1}(y)\frac{d}{dy} = -\half\Phi(\h)		\tag{7.14}
$$
therefore determines $h(y)$ up to left multiplication by 
$h_o\in \tlG = G_{\R}$.  

\par Define
$$
	h(y) = g(y)y^{-H/2}				\tag{7.15}
$$
where $H = \Phi_0(\h)$.  Then, a standard weight argument shows that
$$
	g^{-1}(y)\frac{dg}{dy} 
	= -\half y^{-H/2}(\Phi(\h) - \Phi_0(\h) y^{-1})	
	= \sum_{m\geq 2}\, B_m y^{-2}			\tag{7.16}
$$
Consequently, $g(y)$ and $g^{-1}(y)$ have convergent series expansions
about $\infty$ of the form
$$
\aligned
	g(y) &= g(\infty)(1 + g_1 y^{-1} + g_2 y^{-2} + \cdots)		\\
	g^{-1}(y) &= (1 + f_1 y^{-1} + f_2 y^{-2} + \cdots)g^{-1}(\infty)
\endaligned							\tag{7.17}
$$
where the coefficients $g_k$ and $f_k$ are universal non-commutative
polynomials in the $B_k$ with rational coefficients.  

\par To connect these results with the $SL_2$-orbit theorem, we now
assumes that $\theta(z) = e^{zN}.F$ is a nilpotent orbit of pure
Hodge structure and let
$$
	(F,\rel W) = (e^{-i\delta}.\hat F,\rel W)		\tag{7.18}
$$
be the splitting of the limiting mixed Hodge structure of $\theta(z)$
defined by Theorem $(2.16)$.  Define 
$$
	F_o = \hat\theta(i) = e^{iN}.\hat F
$$
where $\hat\theta(z) = e^{zN}.\hat F$ is the associated split orbit,
and require $\Phi_0$ to be the associated representation of $sl_2(\R)$
defined by Theorem $(3.13)$.  Then,
$$
	h(y).F_o = g(y)y^{-H/2}.F_o = g(y)e^{iyN}.\hat F	
$$
On the other hand, by Theorem $(7.1)$, $h(y).F_o = e^{iy\tilde N}.\tilde F$
and hence 
$$
	e^{iy\tilde N}.\tilde F = g(y)e^{iyN}.\hat F
$$
Therefore, in order to complete the proof of the $SL_2$ orbit theorem, it 
remains only to show that one can select data $(g(\infty),\{T_n\})$ such that
$(\tilde N,\tilde F) = (N,F)$. Assuming that $g(\infty)\in\ker(N)$, this
then boils down after a lengthy calculation to requirement that
$$ 
	e^{i\delta}  = g(\infty)\left(
	      1+\sum_{k>0}\, \frac{1}{k!}(-i)^k(\ad\, N_0)^k\,g_k\right)
$$
At this point, the algebra/combinatorics of solving for $g(\infty)$ and
$\{T_n\}$ becomes sufficiently involved that I shall leave the details to
\S 8 and [CKS].

\head \S 8.\quad Nilpotent Orbits of Type (I) \endhead

\par In this section we prove Theorem $(4.2)$ for admissible nilpotent
orbits of type $(\text{I})$ by constructing a suitable series solution 
$\b(y)$ of the Lax equation $(7.2)$ using the outline of [CKS] developed
in \S 7.  To determine what form the series expansion of $\b(y)$ should
assume, consider the following two examples:

\example{Example 8.1} Let $\pi:E\to\C$ denote the family of elliptic curves
defined by the equation 
$$
	v^2 = u(u-1)(u-s)
$$
and $\tilde\pi:\tilde E\to\C$ denote the corresponding family of punctured
curves obtained by deleting the points of $E$ lying over $u=a$ for some
fixed parameter $a\in\C-\{0,1\}$.  Then, after a local rescaling of 
coordinates, the function
$$
	\b(y) = \Ad(h^{-1}(y))N
$$
attached by Theorem $(6.11)$ to the nilpotent orbit of 
$R^1_{\tilde\pi*}(\Q)\otimes{\Cal O}_{\C-\{0,1,a\}}$ at $s=0$ 
is given by the formula
$$	
	\b(y) = \frac{N}{y} - \frac{\delta}{y^{3/2}}
$$
\endexample

\example{Example 8.2} Let $\hat\theta(z) = e^{zN}.\hat F$ be a split orbit 
of type $(\text{I})$ and $\UU = H(1)\otimes S(1)$ [cf\. Theorem 
$(3.14)$].  Equip $\gg_{\C}$ with the associated $sl_2$-module structure 
defined by Theorem $(3.16)$ and suppose that 
$$
	\Psi:\UU\to\gg_{\C}
$$
is a morphism of Hodge structure with respect to $F_o = \hat\theta(i)$
such that
$$
	\varsigma.\Psi(\tau) = \Psi(\varsigma.\tau)
$$
for all $\varsigma\in sl_2(\C)$ and $\tau\in\frak U$.  Then, 
$$
	\theta(z) = e^{zN}e^{-i\Psi(f)}.\hat F
$$
is an admissible nilpotent orbit of type $(\text{I})$ with split orbit 
$\hat\theta(z)$ and associated functions 
$$
	\b(y) = \frac{N}{y} + \frac{\Psi(f)}{y^{3/2}},\qquad
	h(y) = (1 + \Psi(e)y^{-1})y^{-H/2}
$$
\endexample

\par Accordingly, let us assume that the desired function $\b(y)$
is horizontal with respect to $F_o$ and has a convergent 
series expansion about $\infty$ of the form
$$
	\b(y) = \sum_{n\geq 0}\, \b_n y^{-1-n/2}	\tag{8.3}
$$
Let $\Phi(y)$ be the corresponding function defined by equations
$(6.43)$--$(6.44)$ and $\Psi(y)$ be the linear map from
$\UU = H(1)\otimes S(1)$ to $\gg_{\C}$ defined by the equation
$$
	\Psi(e + if) = 2i\b^{0,-1}(y),\qquad \Psi(e-if) = -2i\b^{-1,0}
							\tag{8.4}
$$
Then, a short calculation shows that equation $(7.2)$ is equivalent
to the pair of differential equations
$$
  -8\Phi'(y) = Q(\Phi,\Phi),\qquad -2\Psi'(y) = Q(\Phi,\Psi)    \tag{8.5}
$$

\par Thus, as in [CKS], the series expansion 
$$
	\Phi(y) = \sum_{n\geq 0}\, \Phi_n y^{-1-n/2}	
$$
of $\Phi$ can be computed inductively starting from a collection of 
morphisms of Hodge structure
$$
	T_n:sl_2(\C)\to\gg_{\C}(n)				\tag{8.6}
$$
such that $\Omega T_n = (n^2-2n) T_n$, where 
$$
	\gg_{\C} = \bigoplus_r\, \gg(r)				\tag{8.7}
$$
denotes the decomposition of $\gg_{\C}$ into isotypical components with
respect to $sl_2$-module structure 
$$
	x.y = [\Phi_0(x),y]
$$
induced by $\Phi_0$ on $\gg_{\C}$.  Moreover, $\Phi_1 = 0$.
\vskip 3pt

\par Similarly, the coefficients of the series expansion 
$$
	\Psi = \sum_{n\geq 0}\, \Psi_n y^{-1-n/2}		
$$
satisfy the recursion relation
$$
	(n+2)\Psi_n = \sum_{j=0}^n\, Q(\Phi_j,\Psi_{n-j})	\tag{8.9}
$$
Therefore, except for the contribution introduced by the term 
$Q(\Phi_0,\Psi_n)$, equation $(8.9)$ allows us to inductively
compute the coefficients of $\Psi$.  
\vskip 3pt

\par To rectify this problem, let $R$ be the endomorphism of 
$\Hom(\UU,\gg_{\C})$ defined by $Q(\Phi_0,*)$ and
recall that if $U_r$ and $U_s$ are irreducible $sl_2$-modules of
 highest weight $r$ and $s$ then 
$$
	U_r\otimes U_s 
	= \bigoplus_{|r-s|<t<r+s,\, t\equiv r+s \mod 2}\, U_t	\tag{8.10}
$$
where $U_t$ is irreducible of highest weight $t$.  In particular,
$$
	\Hom(\UU,\gg_{\C})
	= \Hom(\UU,\gg_{\C})^+\oplus
	  \Hom(\UU,\gg_{\C})^-			\tag{8.11}
$$
where $\Hom(\UU,\gg_{\C})^{\pm}$ is of highest weight $n$
with respect to the left action of $sl_2$ on 
$\Hom(\UU,\gg_{\C}) \cong \gg_{\C}\otimes(\frak U)^*$ and
highest weight $n\pm 1$ with respect to the diagonal action.  

\proclaim{Calculation 8.12} $R$ acts semisimply on 
$\Hom(\UU,\gg(n))$ as multiplication by $(n+2)$ on 
$\Hom(\UU,\gg(n))^-$ and multiplication by $-n$ on 
$\Hom(\UU,\gg(n))^+$. 
\endproclaim
\demo{Proof} Let $e=(1,0)$ and $f=(0,1)$ denote the standard basis of 
$\C^2$ and $M$ be an irreducible submodule of $\gg_{\C}$ of highest
weight $n$.  Then, relative to the standard identification of 
$M$ with $\text{Sym}^n(\C^2)$,
$$
	M\otimes\UU^* \cong A\oplus B			\tag{8.13}
$$
where
$$
\aligned
	A &= \text{span}(a_0,\dots,a_{n+1}),\qquad
	a_j = (n-j+1)e^{n-j}f^j\otimes f^* - j e^{n-j+1}f^{j-1}\otimes e^* \\
	B &= \text{span}(b_0,\dots,b_{n-1}),\qquad
	b_j = e^{n-j-1}f^{j+1}\otimes f^* + e^{n-j}f^j\otimes e^*
\endaligned
$$
are irreducible submodules of highest weight $n+1$ and $n-1$ with respect
to the diagonal action of $sl_2(\C)$, and [cf\. $(3.15)$]
$$
	\h.(a_j) = (n+1-2j)a_j,\qquad
	\h.(b_j) = (n-1-2j)b_j					\tag{8.14}
$$
Accordingly, it suffices to compute $R(a_j)$ and $R(b_j)$.  A short 
calculation shows that
$$
	Q(\sigma,\tau)(v) =  2[\sigma(\n_0^+),\tau(\n_0.v)]
			    +2[\sigma(\n_0^-),\tau(\n_0^+.v)]
			     -[\sigma(\h),\tau(\h.v)]	  
$$
Therefore, 
$$
\aligned
	R(a_j)(e) &= 2\n_0^+.\a_j(f) + \h.\a_j(e)
		   = 2\n_0^+.((n-j+1)e^{n-j}f^j) 
		     +\h.(-je^{n-j+1}f^{j-1})      \\
		  &= 2(n-j+1)je^{n-j+1}f^{j-1} - j(n-2j+2)e^{n-j+1}f^{j-1} \\
		  &= j(2n-2j+2-n+2j-2)e^{n-j+1}f^{j-1}			    
		   = jn e^{n-j+1}f^{j-1}				    
		   = -n a_j(e)						    
\endaligned
$$
The remaining calculation of $R(a_j)(f)$, $R(b_j)(e)$ and $R(b_j)(f)$ are
similar and left to the reader.
\enddemo

\proclaim{Corollary 8.15} $\Psi_0 = 0$, $\Psi_1\in\Hom(\UU,\gg(1))^-$, 
$\Psi_2\in\Hom(\UU,\gg(2))^-$.
\endproclaim
\demo{Proof} By equation $(8.9)$, $R(\Psi_0) = 2\Psi_0$, and hence 
$\Psi_0\in\Hom(\UU,\gg(0))^-$ by Calculation $(8.15)$.  However,
$\Hom(\UU,\gg(0))^- = 0$ since it is highest weight $-1$ with
respect to the diagonal action of $sl_2$.  Consequently, by 
virtue of the fact that $\Psi_0 = 0$ and $\Phi_1 = 0$, it then
follows from equation $(8.9)$ that $R(\Psi_1) = 3\Psi_1$ and
$R(\Psi_2) = 4\Psi_2$.  Therefore, by Calculation $(8.15)$,
$\Psi_1\in\Hom(\UU,\gg(1))^-$ and $\Psi_2\in\Hom(\UU,\gg(2))^-$.
\enddemo

\par To continue, given a semisimple endomorphism of $A$ of a finite
dimensional vector space $V$, let $[*]^A_{\l}$ denote projection from
$V$ onto the $\l$ eigenspace of $V$.  Then, by virtue of Calculation
$(8.15)$,
$$
\aligned
	(n-k) \Psi_{n,k}^- 
	&= \left[\sum_{0<j<n}\, Q(\Phi_j,\Psi_{n-j})\right]^R_{k+2}	\\
	(n+k+2) \Psi_{n,k}^+ 
	&= \left[\sum_{0<j<n}\, Q(\Phi_j,\Psi_{n-j})\right]^R_{-k}
\endaligned\tag{8.16}
$$
where $\Psi_{n,k}^{\pm}$ denotes the component of $\Phi_n$ which takes
values in $\Hom(\UU,\gg(k))^{\pm}$.  Therefore, subject to the compatibility 
condition
$$
  \left[\sum_{0<j<n}\, Q(\Phi_j,\Psi_{n-j})\right]^R_{n+2} = 0	\tag{8.17}
$$
equation $(8.9)$ allows one to compute $\Psi_n$ modulo $\Psi_{n,n}^-$
from $\Phi$ and $\Psi_1,\dots,\Psi_{n-1}$.

\par To handle the compatibility condition $(8.17)$, observe that by 
virtue of equation $(8.10)$,
$$
	\Hom(sl_2(\C),\gg(n)) 
	= \bigoplus_{\epsilon=-1}^1\, \Hom(sl_2(\C),\gg(n))^{\epsilon}
$$
where $\Hom(sl_2(\C),\gg(n))^{\epsilon}$ is highest weight $n$ with 
respect to the left action of $sl_2(\C)$ on $\gg_{\C}$ and highest
weight $n+2\epsilon$ with respect to the diagonal action.

\proclaim{Lemma 8.18} Let $C\in\Hom(sl_2(\C),\gg(r))^{-1}$ and 
$B\in\Hom(\UU,\gg(s))^-$.  Then,
$$
	Q(C,B) \in 
	\bigoplus_{|r-s|\leq t\leq r+s-2,\, t\equiv r+s \mod 2}\, 
	\Hom(\UU,\gg(t))
$$
\endproclaim
\demo{Proof} By equation $(8.10)$ and the Jacobi identity,
$$
	\Span([\gg(r),\gg(s)])
	\subseteq \bigoplus_{|r-s|\leq t\leq r+s,\,t=r+s\mod 2}\, \gg(t)
$$
Therefore, it suffices to show that $Q(C,B)$ projects trivially
onto $\Hom(\UU,\gg(r+s))$. Direct calculation
shows that every irreducible submodule of $\Hom(sl_2(\C),\gg(r))^{-1}$
is isomorphic to $\Span(c_0,...c_{r-2})$ where
$$
	c_k(\n_0)   = e^{r-k-2}f^{k+2},\quad
	c_k(\h) = 2e^{r-k-1}f^{k+1},\quad
	c_k(\n_0^+) = - e^{r-k}f^k
$$
Accordingly, by the semisimplicity of $sl_2(\C)$, it is sufficient to show 
that  
$$
	Q(c_k,b_j)  = 0 \mod \bigoplus_{t\leq r+s-2}\, \Hom(\UU,\gg(t))
$$
Consider $Q(c_k,b_j)(e)$:
$$
\aligned
	Q(c_k,b_j)(e) 
	&= 2[c_k(\n_0^+),b_j(f)] + [c_k(\h),b_j(e)]			\\
	&= -2[e^{r-k}f^k,e^{s-j-1}f^{j+1}] +2[e^{r-k-1}f^{k+1},e^{s-j}f^j] \\
	&\in E_{r+s-2k-2j-2}(\h) 
\endaligned							\tag{8.19}
$$
Suppose that $Q(c_k,b_j)(e)$ projects non-trivially onto $\gg(r+s)$.
Then, by $(8.19)$, $\n_0^{r+s-j-k-1}.Q(c_k,b_j)(e)\neq 0$.  But,
$$
\aligned
	\n_0^{r+s-j-k-1}.&[e^{r-k}f^k, e^{s-j-1}f^{j+1}]		\\
	&= \pmatrix r+s-j-k-1 \\ r-k\endpmatrix
	 [\n_0^{r-k}.e^{r-k}f^k,\n_0^{s-j-1}.e^{s-j-1}f^{j+1}]		\\
	&= \pmatrix r+s-j-k-1 \\ r-k\endpmatrix
	   [(r-k)!f^r,(s-j-1)!f^s]					\\
	&=(r+s-j-k-1)![f^r,f^s]						
\endaligned
$$
Likewise,
$\n_0^{r+s-j-k-1}.[e^{r-k-1}f^{k+1},e^{s-j}f^j] = (r+s-j-k-1)![f^r,f^s]$.
Combining these two equations with $(8.19)$, it then follows that 
$Q(c_k,b_j)(e)$ projects trivially onto $\gg(r+s)$.  Similarly,
one finds that $Q(c_k,b_j)(f)$ projects trivially onto $\gg(r+s)$,
thereby proving the lemma.
\enddemo 

\proclaim{Theorem 8.20} For any choice of a collection of morphisms of Hodge 
structure
$$
	S_n\in\Hom(\UU,\gg(n))^-,\qquad n>0
$$
there exists a unique, convergent $\hh$-valued series solution 
$\Psi = \sum_{n>0}\, \Psi_n y^{-1-n/2}$ of equation $(8.5)$ which is 
horizontal with respect to $F_o$ such that 
\roster
\item"(a)" $\Psi_n\in\oplus_{r\leq n,\, r\equiv n \mod 2}\,\Hom(\UU,\gg(r))$;
\item"(b)" $\Psi_{n,n}^- = S_n$;
\item"(c)" $\Psi_{n,n}^+ = 0$.
\endroster
\endproclaim
\demo{Proof} The desired function $\Psi$ can now be constructed inductively
using equation $(8.16)$.  Namely, by Corollary $(8.15)$, we can assume by 
induction that $\Psi_m$ satisfies conditions $(a)$--$(c)$ for $m<n$.  
Therefore, by Lemma $(8.18)$,
$$
	\sum_{0<j<n}\, Q(\Phi_j,\Psi_{n-j}) \in 
	\bigoplus_{t<n,\, t\equiv n \mod 2}\,\Hom(\UU,\gg(t))	\tag{8.21}
$$
since 
$$
	\Phi_k\in\bigoplus_{s\leq k,\, s\equiv k \mod 2}\,\Hom(sl_2,\gg(s))
$$
by [CKS:6.17].  Consequently, $\sum_{0<j<n}\, Q(\Phi_j,\Psi_{n-j})$ 
satisfies the compatibility condition $(8.17)$, and hence we can
solve for $\Psi_n$ modulo $\Psi_{n,n}^-$ using equation $(8.16)$.
In particular, by equation $(8.21)$ and $(8.16)$, $\Psi_{n,n}^+ =0$.  
Likewise, $\Psi_{n,k} = 0$ for $n>k$.  Thus, given $\Phi$ and 
$\Psi_1,\dots,\Psi_{n-1}$ there exists a unique solution $\Psi_n$ 
to equations $(8.9)$ which satisfies conditions $(a)$--$(c)$.

\par Imposing the condition that $S_n = \Psi_{n,n}^-$ be a morphism
of Hodge structure, it then follows from [CKS:6.47] and equation $(8.16)$
that $\Psi_n$ is horizontal and takes values in $\hh$.

\par To prove that the formal series solution
$$
	\Psi(y) = \sum_{n\geq 0}\, \Psi_n y^{-1-n/2}		\tag{8.21}
$$
constructed above converges about $y=\infty$, recall that $\gg_{\C}$ is a
subalgebra of $gl(V)$ and let $||*||$ be norm on $gl(V)$ such that
$||AB|| \leq ||A|| ||B||$.  Define
$$
\aligned
	||A||_1 &= 4(||A(\xx^+)|| + ||A(\xx^-)|| + ||A(\zz)||)	\\
	||B||_2 &= ||B(\nu_+)|| + ||B(\nu_-)||
\endaligned\qquad
\aligned
	A &\in\Hom(sl_2(\C),\gg_{\C})	\\
	B &\in\Hom(\UU,\gg_{\C})
\endaligned	
$$
Then, a short calculation shows that
$$
	||Q(A,B)||_2 \leq ||A||_1 ||B||_2
$$
Therefore, by equation $(8.9)$, 
$$
	(n+2)||\Psi_n||_2 \leq ||\Phi_0||_1 ||\Psi_n||_2 
	                     + \sum_{0<j<n}\, ||\Phi_j||_1 ||\Psi_{n-j}||_2
$$
and hence
$$
	(n-1)||\Psi_n||_2 \leq \sum_{0<j<n}\, ||\Phi_j||_1 ||\Psi_{n-j}||_2
								\tag{8.22}
$$
upon rescaling $||*||$ so that $||\Phi_0||_1 = 3$.  

\par To continue, we note that since $\gg_{\C}$ is finite dimensional,
there exists an integer $m$ such that $\gg(n) = 0$ for $n>m$.  
Consequently, $S_n = 0$ for $n>m$ and hence
$$
	\max_k\, ||S_k||_2
$$
is finite.  Therefore, there exists a constant $D$ such that\footnote
{In the degenerate case $\max_k\, ||S_k||_2 =0$ all $S_k = 0$ and
hence $\Psi = 0$ by Theorem $(8.20)$.} 
$$
    ||\Psi_{\ell}||_2 \leq D^{\ell}(\max_k\, ||S_k||_2)^{\ell}	\tag{8.23}
$$
for $\ell\leq m$.  Similarly, by [CKS:6.24] there exists a constant
$C$ such that 
$$
    ||\Phi_{\ell}||_1 \leq  C^{\ell}(\max_k ||T_k||_1)^{\ell}
$$

\par Assume by induction that $(8.23)$ holds for $\ell<n$, and enlarge
$D$ if necessary so that 
$$
	D(\max_k\, ||S_k||_2) \geq C(\max_k ||T_k||_1)
$$
Then, by equation $(8.22)$,
$$
\aligned	
	(n-1)||\Psi_n||_2 
	&\leq \sum_{0<j<n}\, ||\Phi_j||_1 ||\Psi_{n-j}||_2	\\
	&\leq \sum_{0<j<n}\, C^j(\max_k ||T_k||_1)^j 
			     D^{n-j}(\max_k ||S_k||_2)^{n-j}    \\
        &\leq \sum_{0<j<n}\, D^n(\max_k ||S_k||_2)^n		
	 = (n-1)D^n(\max_k ||S_k||_2)^n
\endaligned
$$
Therefore,
$$
	||\Psi_n||_2 \leq D^n(\max_k ||S_k||_2)^n
$$ 
for all $n$, and hence the series $(8.21)$ converges on some interval 
$(a,\infty)$.
\enddemo

\par Invoking Theorem $(7.1)$, we now obtain an $\tlG$-valued function
$h(y)$ such that
$$
    h^{-1}\frac{dh}{dy} 
    = -\LL\b(y) 
    = -\half\Phi(\h) - \Psi(e)					\tag{8.24}
$$
Following [CKS], let $H = \Phi_0(\h)$ and $g(y)$ be the $\tlG$-valued
function defined by the equation
$$
	h(y) = g(y)y^{-H/2}					\tag{8.25}
$$
Then,
$$
\aligned
	\left[g^{-1}\frac{dg}{dy}\right]^{\ad\, Y}_0
	&= -\half y^{-H/2}.(\Phi(\h) - \Phi_0(\h)y^{-1})	\\	
	\left[g^{-1}\frac{dg}{dy}\right]^{\ad\, Y}_{-1}
	&= -y^{-H/2}.\Psi(e)					
\endaligned							\tag{8.26}
$$	
where $Y=Y_{(F_o,W)}$.

\proclaim{Theorem 8.27} $g^{-1}(dg/dy) = \sum_{m\geq 2}\, B_m y^{-m}$.
\endproclaim
\demo{Proof} Due to the short length of $W$,
$$
	g^{-1}\frac{dg}{dy} 
	= \left[g^{-1}\frac{dg}{dy}\right]^{\ad\, Y}_0 + 
	  \left[g^{-1}\frac{dg}{dy}\right]^{\ad\, Y}_{-1}
$$
Therefore, since $\Phi$ is isomorphic via the grading $Y$ with the 
corresponding function defined by nilpotent orbits of pure Hodge
structure induced by $\theta(z)$ on $Gr^W$, it then follows from 
[CKS:6.30] that
$$
	\left[g^{-1}\frac{dg}{dy}\right]^{\ad\, Y}_0	
	= \sum_{m\geq 2}\, [B_m]^{\ad\,Y}_0 y^{-m}
$$
where
$$
	\left[B_m\right]^{\ad\,Y}_0 
	= -\half\sum_{n\geq m}\, [\Phi_n(\h)]^{\ad\,H}_{2(m-1)-n}
$$

\par To establish that $[g^{-1}\frac{dg}{dy}]^Y_{-1}$ is
also of this form, observe that by $(8.26)$:
$$
\aligned
	\left[g^{-1}\frac{dg}{dy}\right]^Y_{-1} 
	&= -y^{-H/2}.\Psi(e) 
	 = -y^{-H/2}.\left(\sum_{n>0}\, \Psi_n(e)\, y^{-1-n/2}\right)	  \\
	&= -y^{-H/2}.\left(\sum_{n>0}\sum_{r=0}^n\, 
		[\Psi_n(e)]^H_{n-2r} y^{-1-n/2}\right)			  \\
	&= -\sum_{n>0}\sum_{r=0}^n\, [\Psi_n(e)]^H_{n-2r}\, y^{-1-n+r} 
\endaligned							\tag{8.28}
$$
However, by the description of the irreducible submodules $B$ of 
$\Hom(\UU,\gg(n))^-$ presented in Calculation $(8.12)$,
$[\Psi_n(e)]^H_{-n} = 0$ and hence equation $(8.28)$ reduces to 
$$
	\left[g^{-1}\frac{dg}{dy}\right]^Y_{-1} 
	= -\sum_{n>0}\sum_{r=0}^{n-1}\, [\Psi_n(e)]^H_{n-2r}\, y^{-1-n+r}
	=\sum_{m\geq 2}\, [B_m]^{\ad\, Y}_{-1} y^{-m}
$$
where
$$
  	[B_m]^{\ad\,Y}_{-1} = 
 	 -\sum_{n\geq m-1}\, [\Psi_n(e)]^{\ad\,H}_{2(m-1)-n}	\tag{8.29}
$$
\enddemo	

\proclaim{Corollary 8.30} The functions $g(y)$ and $g^{-1}(y)$ have 
convergent Taylor expansions about $y=\infty$ of the form 
$$
\aligned
	g(y) &= g(\infty)(1 + g_1 y^{-1} + g_2 y^{-2} + \cdots)		\\
	g^{-1}(y) &= (1 + f_1 y^{-1} + f_2 y^{-2} + \cdots)g^{-1}(\infty)
\endaligned
$$
where $g(\infty)$ is an arbitrary element of $\tlG$ determined by the
initial value of $h(y)$.  Moreover, the coefficients $g_n$ and $f_n$ can be 
expressed as universal non-commutative polynomials in the $B_k$ with rational 
coefficients, weighted homogeneous of degree $n$ when $B_k$ when $B_k$ is 
assigned weight $k-1$.  $B_{n+1}$ occurs with coefficient $-1/n$ in $g_n$ 
and with coefficient $1/n$ in the case of $f_n$.
\endproclaim
\demo{Proof} See Lemma $(6.32)$ in [CKS].
\enddemo 

\proclaim{Calculation 8.31}  $\n_0^k. B_k = 0$.
\endproclaim
\demo{Proof} That $\n_0^k.[B_k]^Y_0 = 0$ is shown in [CKS:6.32].
Moreover, by $(8.29)$:
$$
	[B_k]^Y_{-1} = -\sum_{n\geq k-1}\, [\Psi_n(e)]^H_{2(k-1)-n}
$$
and hence
$$
	\n_0^k. [B_k]^Y_{-1} 
	= -\sum_{n\geq k-1}\, \n_0^k. [\Psi_n(e)]^H_{2(k-1)-n} = 0
$$
since $\Psi_n(e)$ takes values in $\oplus_{r\leq n}\, \gg(r)$.
\enddemo

\proclaim{Corollary 8.32}  
$\n_0^{k+1}. g_k = \n_0^{k+1}.f_k = 0$.
\endproclaim 
\demo{Proof} By Corollary $(8.30)$, $g_k$ and $f_k$ are homogeneous polynomials
of degree $k$ in $B_2,\dots, B_{k+1}$ with respect to the grading  
$\deg(B_{\ell}) = \ell - 1$.  Therefore, by virtue of Calculation $(8.31)$ and
Leibniz rule, both $\n_0^{k+1}. g_k$ and $\n_0^{k+1}. f_k = 0$.
\enddemo

\proclaim{Theorem 8.33} Let $\b(y) = \Phi(\n_0) + \Psi(f)$ denote the 
solution equation $(7.2)$ constructed above, and 
$e^{z\tilde N}.\tilde F$ be the associated nilpotent orbit defined by 
Theorem $(7.1)$.  Then, $\tilde N$ coincides with $N_0 = \Phi_0(\n_0)$ if 
and only if $g(\infty)\in\ker(\ad\, N_0)$.
\endproclaim
\demo{Proof:} By definition, 
$$
	    \tilde N = h(y).\b(y) 
		     = h(y).([\b(y)]^Y_0 + [\b(y)]^Y_{-1})
		     = h(y).[\b(y)]^Y_0 + h(y).\psi(f)
$$ 
Moreover, since $\Psi_0 = 0$, 
$$
\aligned
	y^{-H/2}.\Psi(f) 
	&= y^{-H/2}.\left(\sum_{n>0}\,\sum_{r=0}^n\, 
		[\Psi_n(f)]^H_{n-2r}\, y^{-1-n/2}\right)	\\
	&= \sum_{n>0}\,\sum_{r=0}^n\,
                [\Psi_n(f)]^H_{n-2r}\, y^{-1-n+r}		
	= \{\cdots\}y^{-1} + \{\cdots\}y^{-2} + \cdots
\endaligned
$$
Thus, making use of the calculations of [CKS], we have
$$
\aligned
    \tilde N &=  h(y).[\b(y)]^Y_0 + \{\cdots\}y^{-1} + \{\cdots\}y^{-2} 
				  + \cdots				   \\
	     &=  g(y)y^{-H/2}.[\b(y)]^Y_0 + \{\cdots\}y^{-1} 
		 + \{\cdots\}y^{-2} + \cdots				   \\
	     &=  g(y).(N_0 + \{\cdots\}y^{-1} + \{\cdots\}y^{-2} + \cdots)
		 + \{\cdots\}y^{-1} + \{\cdots\}y^{-2} + \cdots		   \\
	     &=  g(\infty).N_0 + \{\cdots\}y^{-1} + \{\cdots\}y^{-2} 
		 + \cdots
\endaligned
$$
and hence $\tilde N = g(\infty).N_0$.
\enddemo


\par To connect previous constructions with Theorem $(4.2)$, let us now
suppose that $\theta(z) = e^{zN}.F$ is an admissible nilpotent orbit of 
type $(\text{I})$, and let 
$$
	\hat\theta(z) = e^{zN}.\hat F
$$
be the associated split orbit obtained by applying the splitting operation
$$
	(\hat F,\rel W) = (e^{-i\delta}.F,\rel W)
$$
to the limiting mixed Hodge structure of $\theta$.  Define
$$
	F_o = \hat\theta(\sqrt{-1})				\tag{8.34}
$$
and let $(N_0,H,N_0^+)$ be the associated $sl_2$-triple obtained by 
application of Theorem $(3.16)$ to $\hat\theta$.  Set
$$
	\Phi_0(\n_0) = N_0,\qquad \Phi(\h) = H,\qquad
	\Phi_0(\n_0^+) = N_0^+					\tag{8.35}
$$
and recall that $N_0 = N$ due to the short length of $W$.

\proclaim{Theorem 8.36} Let $\b(y) = \Phi(\n_0) + \Psi(f)$ denote the 
solution equation $(7.2)$ constructed above, and $e^{z\tilde N}.\tilde F$ 
be the associated nilpotent orbit obtained from Theorem $(7.1)$.  Assume
that $F_o$ and $\Phi_0$ are given by equations $(8.34)$--$(8.35)$ and
$g(\infty)\in\ker(\ad\,N_0)$.  Then,
$$
	\tilde F 
	= g(\infty)\left(
	      1+\sum_{k>0}\, \frac{1}{k!}(-i)^k(\ad\, N_0)^k\,g_k\right).
	      \hat F
$$
\endproclaim
\demo{Proof}  By Theorem $(7.1)$, $h(y).F_0 = e^{iyN_0}.\tilde F$.
Therefore, 
$$
\aligned
   \tilde F &= e^{-iyN_0}h(y) 
	     = e^{-iyN_0}g(\infty)
	       \left(1 + \sum_{k>0}\, g_k y^{-k}\right)
	       y^{-H/2}e^{iN_0}.\hat F					\\
	    &= e^{-iyN_0}g(\infty)
               \left(1 + \sum_{k>0}\, g_k y^{-k}\right)
               e^{iyN_0}.\hat F						\\
	    &= g(\infty)e^{-iyN_0}
               \left(1 + \sum_{k>0}\, g_k y^{-k}\right)
               e^{iyN_0}.\hat F						\\
 	    &= g(\infty)\left(e^{-iy\, \ad\,N_0}\,
               \left(1 + \sum_{k>0}\, g_k y^{-k}\right)\right)
		.\hat F		  					\\
	    &= g(\infty)
               \left(1 + \sum_{k>0,j\geq 0}\, 
			    \frac{(-i)^j}{j!}(ad N_0)^j 
			     g_k y^{j-k}\right).\hat F	
\endaligned
$$
Moreover, by Corollary $(8.32)$, $(ad\, N_0)^j\, g_k = 0$ whenever $j>k$.
Thus,
$$
\aligned
   \tilde F &= g(\infty)
               \left(1 + \sum_{k>0}\sum_{j=0}^k\,
                            \frac{(-i)^j}{j!}(ad N_0)^j
                             g_k y^{j-k}\right).\hat F_{\infty}  	\\
	    &= g(\infty)\left(1+\sum_{k>0}\,
                 \frac{1}{k!} (-i)^k (ad\,N_0)^k\, g_k\right).\hat F_{\infty}
		 + \{\cdots\}y^{-1} + \{\cdots\}y^{-2} + \cdots
\endaligned 
$$
Accordingly, upon taking the limit as $y\to\infty$ in this last equation
we obtain the stated formula for $\tilde F$.
\enddemo

\par Thus, in order to complete the proof of Theorem $(4.2)$ for admissible
nilpotent orbits of type $(\text{I})$, it is sufficient to show that 
we can select morphisms of Hodge structure
$$
	T_n\in\Hom(sl_2(\C),\gg(n))^{-1},\qquad
	S_n\in\Hom(\UU,\gg(n))^-
$$
for $n>0$ and element 
$\zeta = \log(g(\infty))\in\hh\cap\ker(\ad\,N_0)\cap\lam_{(\hat F,\rel W)}$ 
such that
$$ 
	e^{i\delta}  = e^{\zeta}\left(
	      1+\sum_{k>0}\, \frac{1}{k!}(-i)^k(\ad\, N_0)^k\,g_k\right)
$$

\proclaim{Theorem 8.37} Let $\theta(z) = e^{zN}.F$ be an admissible 
nilpotent orbit of type $(\text{I})$.  Then, the solutions $\b(y)$
of equation $(7.2)$ which have the following three properties
\roster
\item $\b(y)$ is horizontal at $F_o = \hat\theta(i)$;
\item $\b(y) = \sum_{n\geq 0}\, \b_n y^{-1-n/2}$;
\item $\b_0 = N_0$;
\endroster
are in 1-1 correspondence with the elements
$
	\eta \in \hh\cap\ker(\ad\, N_0)\cap\lam_{(\hat F,\rel W)}
$
via the map
$$
	\eta = \sum_{n>0}\, [\b_n]^{\ad\, H}_{-n}
$$
\endproclaim
\demo{Proof} If $\b(y)$ satisfies the conditions stated above then so
does $[\b(y)]^{\ad\, Y}_0$.  Therefore, by Lemma $(6.41)$ in [CKS] the map
$$
	[\eta]^{\ad\,Y}_0 
	= \left[\sum_{n>0}\, [\b_n]^{\ad\, H}_{-n}\right]^{\ad\, Y}_0 
	= \sum_{n>0}\, \left[[\b_n]^{\ad\,Y}_0\right]^{\ad\, H}_{-n} 
$$
determines a bijective correspondence between the morphisms $T_n$ and 
the elements of 
$
	\hh\cap\ker(\ad\,Y)\cap\ker(\ad\, N_0)\cap\lam_{(\hat F,\rel W)}
$.

\par To recover the morphisms $S_n$ from $[\eta]^{\ad\,Y}_{-1}$, 
observe that since $(F_o,W)$ is split over $\R$,
$$
	{\eusm H} = \bigoplus_{r+s=-1}\,\gg^{r,s}_{(F_o,W)}
$$
is a pure Hodge structure of weight $-1$ with respect to which the 
representation of $sl_2(\C)$ defined by $\ad\,\Phi_0$ is Hodge.
Therefore, by Theorem $(3.14)$ we can decompose ${\eusm H}$ into
a direct sum of irreducible submodules $M$, each of which is isomorphic
to one of the following two standard types
\roster
\item"(a)" $H(d)\otimes S(n)$, $n=2d-1$ odd;
\item"(b)" $E^{p,q}\otimes S(n)$, $n+p+q=-1$, $p-q>0$;
\endroster
where $S(n) = \text{Sym}^n(\C^2)$ is the standard representation of 
$sl_2(\C)$ of highest weight $n$ equipped with the Hodge structured
obtained by declaring 
$$
	\nu_r = (e+if)^r(e-if)^{n-r}				\tag{8.38}
$$
to be of type $(r,n-r)$ and $H(d) = \Span(\epsilon^{-d,-d})$ and 
$E(p,q) = \Span(\epsilon^{p,q},\epsilon^{q,p})$ are trivial 
representations of $sl_2$ equipped with the Hodge structure obtained by
requiring $\epsilon^{r,s}$ to type $(r,s)$ and 
$\overline{\epsilon^{r,s}} = \epsilon^{s,r}$.

\par Let $S_n^M$ denote the projection of $S_n$ onto such an irreducible
module $M$.  Then, a short calculation shows that 
$$
	S_n^M(e+if) = \tau_M\epsilon^{-d,-d}\otimes\nu_d,\qquad
	S_n^M(e-if) = \tau_M\epsilon^{-d,-d}\otimes\nu_{d-1}	\tag{8.39}
$$
for some real number $\tau_M$ if $M$ is of type $(a)$.  Similarly,
if $M$ is type $(b)$ then 
$$
\aligned
	S_n^M(e+if) &= \tau_M\epsilon^{p,q}\otimes\nu_{-p}
		      + \bar \tau_M\epsilon^{q,p}\otimes\nu_{-q}	\\
        S_n^M(e-if) &= \tau_M\epsilon^{p,q}\otimes\nu_{-p-1} 
	              + \bar\tau_M\epsilon^{q,p}\otimes\nu_{-q-1}
\endaligned							\tag{8.40}
$$
where $\tau_M\in\C$, $p,q<0$ and $p+q+n = -1$.

\par In particular, if $S_n^M$ is of type $(8.40)$ then
$$
	2iS^M_n(f) 
	= \tau_M\epsilon^{p,q}\otimes(\nu_{-p}-\nu_{-p-1})
	+\bar\tau_M\epsilon^{q,p}\otimes(\nu_{-q}-\nu_{q-1})	
$$
Moreover, for any index $0\leq k\leq n$, 
$$
\aligned
	\nu_k - \nu_{k-1} 
	&= (e+if)^k(e-if)^{n-k} - (e+if)^{k-1}(e-if)^{n-k+1}		\\
	&= i^k (-i)^{n-k} f^n - i^{k-1}(-i)^{n-k+1}f^n + e(\cdots)	\\
	&= (2i) i^{2k-n-1} f^n + e(\cdots)
\endaligned						
$$
Accordingly, using the identity $p+q+n=-1$, it then follows that
$$
	[\b_n^M]^{\ad\,H}_{-n} 
	= [S^M_n(f)]^{\ad\,H}_{-n}
	= (-i)^{\chi}\tau_M\epsilon^{p,q}\otimes f^n
	  + i^{\chi}\bar\tau_M\epsilon^{q,p}\otimes f^n		\tag{8.41}
$$
where $\chi = p-q$.  Similarly, if $S^M_n$ is of type $(8.39)$ then
$$
    [\b_n^M]^{\ad\,H}_{-n} = \tau_M\epsilon^{-d,-d}\otimes f^n 	\tag{8.42}
$$
Therefore, the sum 
$$
	[\eta]^{\ad\,Y}_{-1} 
	= \sum_M\, \eta^M = \sum_{n>0}\sum_M\, [\b_n^M]^{\ad\,H}_{-n}  
								\tag{8.43}
$$
determines $\tau_M$ for all $M$.  

\par To verify that the sum $(8.43)$ takes values in $\lam_{(\hat F,\rel W)}$,
suppose that $S^M_n$ is of type $(8.40)$ and observe that
$$
	e^{iN_0}.(\epsilon^{p,q}\otimes e^n)	
	= \epsilon^{p,q}\otimes\nu_n \in \gg^{n+p,q}_{(F_o,W)}
$$
and hence 
$$
\aligned
	\{e^{iN_0}.(\epsilon^{p,q}\otimes e^n)\}(F_o^r)	
	&= e^{iN_0}(\epsilon^{p,q}\otimes e^n)e^{-iN_0}e^{iN_0}.\hat F^r \\
	&= e^{iN_0}(\epsilon^{p,q}\otimes e^n).\hat F^r			 
	 \subseteq e^{iN_0}.\hat F^{n+p+r}
\endaligned
$$
Therefore,
$$
   (\epsilon^{p,q}\otimes e^n)(\hat F^r)\subseteq \hat F^{n+p+r}   \tag{8.44}
$$
Furthermore, by Theorem $(3.16)$,
$$
	H = \rel Y - Y
$$
where $\rel Y$ is the grading of $\rel W$ defined by the $I^{p,q}$'s of 
$(\hat F,\rel W)$ and $Y$ is the grading of $W$ defined by the $I^{p,q}$'s
of $(F_o,W)$.  Consequently, the condition that $\epsilon^{p,q}\otimes e^n$
be of weight $n$ with respect to $H$ and weight $-1$ with respect to $Y$
implies that 
$$
	\epsilon^{p,q}\otimes e^n 
	\in \bigoplus_t\, \gg^{t,n-1-t}_{(\hat F,\rel W)}		
$$
Imposing the condition $(8.44)$, it then follows that
$$	
	\epsilon^{p,q}\otimes e^n 
	\in \bigoplus_{t\geq n+p}\, \gg^{t,n-1-t}_{(\hat F,\rel W)}  \tag{8.45}
$$
Likewise, switching the roles of $p$ and $q$,
$$
	\epsilon^{q,p}\otimes e^n 
	\in \bigoplus_{s\geq n+q}\, \gg^{s,n-1-s}_{(\hat F,\rel W)}  \tag{8.46}
$$	
Thus, since $\overline{\epsilon^{q,p}\otimes e^n} = \epsilon^{p,q}\otimes e^n$
and $(\hat F,\rel W)$ is split over $\R$, equations $(8.45)$ and $(8.46)$ 
imply that the Hodge components
$$
	(\epsilon^{p,q}\otimes e^n)^{t,n-1-t} 
$$
of $\epsilon^{p,q}\otimes e^n$ with respect to $(\hat F,\rel W)$ vanish
unless
$$
	t = n-1-s,\qquad t\geq n+p,\qquad s\geq n+q		\tag{8.47}
$$
Recalling that $p+q+n = -1$, it then follows from equation $(8.47)$ that
$$
	(\epsilon^{p,q}\otimes e^n)^{t,n-1-t} = 0
$$ unless $t = n+p$.  Accordingly,  since $N_0$ is a $(-1,-1)$-morphism of 
$(\hat F,\rel W)$, 
$$
	\epsilon^{p,q}\otimes f^n = (N_0)^n.(\epsilon^{p,q}\otimes e^n)
	\in\gg^{p,q}_{(F_o,\rel W)}				\tag{8.48}
$$
Now, by equation $(8.40)$, $p,q<0$.  Therefore, by equation $(8.41)$ 
and $(8.48)$,
$$
	[S_n^M]^{\ad\, H}_{-n}
	\in\gg^{p,q}_{(\hat F,\rel W)}\oplus\gg^{q,p}_{(\hat F,\rel W)}
	\subseteq\lam_{(\hat F,\rel W)}
$$
Similarly, if $S_n^M$ is of type $(8.39)$ then 
$$
	[S_n^M]^{\ad\, H}_{-n}
	\in \gg^{-d,-d}_{(\hat F,\rel W)}
	\subseteq\lam_{(\hat F,\rel W)}
$$

\enddemo

\par Following [CKS], we now note that by virtue of Corollary $(8.30)$
$$ 
	1+\sum_{k>0}\, \frac{1}{k!}(-i)^k(\ad\, N_0)^k\,g_k
	= \exp\left(\sum_{k>0}\, Q_k(C_2,\dots,C_{k+1})\right)	 \tag{8.49}
$$
where
$C_{\ell+1} = \frac{(-i)^{\ell}}{\ell!}(\ad\,N_0)^{\ell} B_{\ell+1}$.

\proclaim{Calculation 8.51} Let $(1-x)^r(1+x)^s = \sum_t\, b_{r,s}^t x^t$.
Then, 
$$
	[C_{\ell+1}]^{\ad\,Y}_0
	= i\sum_{p,q\geq 1,p+q\geq\ell}\,
	  ([\eta]^{\ad\, Y}_0)^{-p,-q}
$$
where $([\eta]^{\ad\, Y}_0)^{-p,-q}$ denotes the component of 
$[\eta]^{\ad\, Y}_0$ of type $(-p,-q)$ with respect to $(\hat F,\rel W)$.
\endproclaim
\demo{Proof} See Lemma $(6.60)$ in [CKS].
\enddemo

\proclaim{Calculation 8.52} $[C_{\ell+1}]^{\ad\,Y}_{-1} 
=  i\sum_{p,q\geq 1, p+q\geq \ell+1}\, 
	  b^{\ell-1}_{p-1,q-1} ([\eta]^{\ad\,Y}_{-1})^{-p,-q}$.
\endproclaim
\demo{Proof} By equation $(8.29)$, 
$$
\aligned
	[C_{\ell+1}]^{\ad\, Y}_{-1} 
	&= -\frac{(-i)^{\ell}}{\ell!}(\ad\,N_0)^{\ell}
	    \sum_{n\geq \ell}\, [\Psi_n(e)]^{\ad\,H}_{2\ell-n}	\\
	&= -\frac{(-i)^{\ell}}{\ell!}
            \sum_{n\geq \ell}\, (\ad N_0)^{\ell}[S_n(e)]^{\ad\,H}_{2\ell-n}  \\
	&= -\frac{(-i)^{\ell}}{\ell!}
            \sum_{n\geq \ell}\sum_M\, 
	    (\ad N_0)^{\ell}[S_n^M(e)]^{\ad\,H}_{2\ell-n}
\endaligned							\tag{8.53}
$$
where $S_n = \sum_M\, S_n^M$ denotes the decomposition of $S_n$ into
irreducible components of type $(8.39)$ and $(8.40)$.  

\par Now, for any index $0\leq k\leq n$, 
$$
\aligned
   \nu_k &= (e+if)^k(e-if)^{n-k}
          = (i(f-ie))^k((-i)(f+ie))^{n-k} 			\\
	 &= i^{2k-n}(f-ie)^k(f+ie)^{n-k}
	  = i^{2k-n}\sum_t\, i^t b_{k,n-k}^t e^t f^{n-t}
\endaligned							\tag{8.54}
$$
Therefore, if $S_n^M$ is of type $(8.40)$ then 
$$
\aligned
      [S_n^M(e)]^{\ad\,H}_{2\ell-n}
      &=  \half\tau_M\epsilon^{p,q}
          \otimes[\nu_{-p}+\nu_{-p-1}]^{\ad\,H}_{2\ell-n}
       +\half\bar\tau_M\epsilon^{q,p}
           \otimes[\nu_{-q}+\nu_{-q-1}]^{\ad\,H}_{2\ell-n} \\
      &=  \half\tau_M\epsilon^{p,q}\otimes
           \left(i^{-2p-n+\ell}b_{-p,n+p}^{\ell} 
	         + i^{-2p-2+n+\ell}b_{-p-1,n+p+1}^{\ell}\right)
		 e^{\ell}f^{n-\ell}					\\
      &+\half\bar\tau_M\epsilon^{q,p}\otimes
	  \left(i^{-2q-n+\ell}b_{-q,n+q}^{\ell} 
	         + i^{-2q-2+n+\ell}b_{-q-1,n+q+1}^{\ell}\right)
		 e^{\ell}f^{n-\ell}					\\
      &= \half i^{1+\ell-\chi}\tau_M \epsilon^{p,q}\otimes
         \left(b_{-p,-q-1}^{\ell} - b_{-p-1,-q}^{\ell}\right)
	 e^{\ell}f^{n-\ell}		\\
      &+ \half i^{1+\ell+\chi}\bar\tau_M\epsilon^{q,p}\otimes
         \left(b_{-q,-p-1}^{\ell} - b_{-q-1,-p}^{\ell}\right)
	 e^{\ell}f^{n-\ell}
\endaligned
$$
where $\chi = p-q$ [recall: p+q+n=-1].  To simplify the above expression,
observe that 
$$
\aligned
     \sum_t\, (b_{k,n-k}^t - b_{k-1,n-k+1}^t) x^t
     &= (1-x)^k (1+x)^{n-k} - (1-x)^{k-1} (1+x)^{n-k+1}		\\
     &= (1-x)^{k-1}(1+x)^{n-k}((1-x) - (1+x))			\\
     &= (-2x) (1-x)^{k-1} (1+x)^{n-k}				\\
     &= (-2x) \sum_t\, b_{k-1,n-k}^t x^t
\endaligned
$$
and hence
$$
\aligned
	b_{-p,-q-1}^{\ell} - b_{-p-1,-q}^{\ell}
	&= -2 b^{\ell-1}_{-p-1,-q-1}				\\
	b_{-q,-p-1}^{\ell} - b_{-q-1,-p}^{\ell}
	&= -2 b^{\ell-1}_{-q-1,-p-1}
\endaligned
$$
Accordingly,
$$
\aligned
	[S_n^M(e)]^{\ad\,H}_{2\ell-n}
	 &= -b^{\ell-1}_{-p-1,-q-1}
	 (i^{1+\ell-\chi}\tau_M \epsilon^{p,q}\otimes e^{\ell}f^{n-\ell})  \\
         &-b^{\ell-1}_{-q-1,-p-1}
	   (i^{1+\ell+\chi} \bar\tau_M
	   \epsilon^{q,p}\otimes e^{\ell}f^{n-\ell})
\endaligned							\tag{8.55}
$$
Inserting $(8.55)$ into equation $(8.53)$ it then follows by equation 
$(8.48)$ that
$$
\aligned
	C^M_{\ell+1} 
	&= i b^{\ell-1}_{-p-1,-q-1} i^{-\chi}\tau_M
	   \epsilon^{p,q}\otimes f^n +
           i b^{\ell-1}_{-p-1,-q-1} i^{\chi}\bar\tau_M
	   \epsilon^{q,p}\otimes f^n				\\
        &=  i b^{\ell-1}_{-p-1,-q-1} (\eta^M)^{p,q} +
            i b^{\ell-1}_{-q-1,-p-1} (\eta^M)^{q,p}
\endaligned							\tag{8.56}
$$

\par Similarly, if $S_n^M$ is of type $(8.39)$ then
$$
	C^M_{\ell+1} = i b^{\ell-1}_{d-1,d-1}(\eta^M)^{-d,-d}   \tag{8.57}
$$
Thus, combining equations $(8.56)$ and $(8.57)$ and switching the signs
of $p$ and $q$, we obtain the formula:
$$
	[C_{\ell+1}]^{\ad\,Y}_{-1} 
	= \sum_M\, C_{\ell+1}^M 
	= i\sum_{p,q\geq 1, p+q\geq \ell+1}\, 
	  b^{\ell-1}_{p-1,q-1} ([\eta]^{\ad\,Y}_{-1})^{-p,-q}
$$
\enddemo

\par In particular, by virtue of Calculations $(8.51)$ and $(8.52)$,
$$
	C_{\ell+1} =  i\sum_{p,q\geq 1, p+q\geq \ell+1}\, 
	  b^{\ell-1}_{p-1,q-1} \eta^{-p,-q}			
$$
Therefore,  since $C_{\ell+1}$ is of the same algebraic form as in 
Lemma $(6.60)$ of [CKS], we can use this result verbatim to prove
that given 
$
	\delta\in\hh\cap\ker(N)\cap\lam_{(\hat F,\rel W)}
$
we can find unique elements 
$
	\zeta,\hph{a}\eta\in \hh\cap\ker(N)\cap\lam_{(\hat F,\rel W)}
$ 
such that 
$$ 
	e^{i\delta}  = e^{\zeta}\left(
	      1+\sum_{k>0}\, \frac{1}{k!}(-i)^k(\ad\, N_0)^k\,g_k\right)
$$
By the above remarks, this completes the proof of Theorem $(4.2)$ for
admissible orbits of type $(\text{I})$.

\head \S 9.\quad Nilpotent Orbits of Type (II) \endhead

\par Suppose now that $\theta(z) = e^{zN}.F$ is an admissible nilpotent
orbit of type $(\text{II})$ and let $\hat\theta(z) = e^{zN}.\hat F$ be
the associated split orbit.  Then, application of Theorem $(3.16)$ to
$\hat\theta(z)$ defines a corresponding splitting 
$$
	N=N_0 + N_{_2}						\tag{9.1}
$$
of $N$ such that $\hat\theta_0(z) = e^{zN_0}.\hat F$  is an 
$\text{\rm SL}_2$-orbit. Consequently, 
$$
	F_o = \hat\theta_0(i) \in\M_{\R}			
$$
Furthermore, since $\theta(z)$ is of type $(\text{II})$ the Hodge
decomposition of the associated function $\b(y) = \Ad(h^{-1}(y))N$
defined by Theorem $(6.11)$ is of the form
$$
	\b(y) = \b^{1,-1} + b^{0,0} + \b^{-1,1} + \b^{0,-1} + \b^{-1,0}
               + \b^{-1,-1}					\tag{9.2}
$$
As in \S 7--8, the first five components of the right hand side of equation
$(9.2)$ are governed by the system of differential equations
$$
	-8\Phi' = Q(\Phi,\Phi),\qquad	-2\Psi' = Q(\Phi,\Psi)
$$
Therefore, as in \S 7--8, we can formally solve for these components 
starting from a collection of morphisms of Hodge structures
$$
	T_n:sl_2(\C)\to\gg_{\C},\qquad S_n:\UU\to\gg_{\C}
$$

\par To solve for $\b^{-1,-1}$, we now return to equation $(6.41)$, which
implies that
$$
	\frac{d}{dy}\b^{-1,-1} 
	=   i[\b^{0,0},\b^{-1,-1}] +2i[\b^{0,-1},\b^{-1,0}]	\tag{9.3}
$$
Next, we recall that since $\theta(z)$ is of type $(\text{II})$ there
exists an index $k$ such that the Hodge decomposition of $(F_o,W)$ is
of the form
$$
	V = I^{k,k}\oplus\left(\bigoplus_{p+q=2k-1}\, I^{p,q}\right)
            \oplus I^{k-1,k-1}					\tag{9.4}
$$
for some index $k$.  Therefore, since $Gr^W_{2k}$ and $Gr^W_{2k-2}$ are
of pure type $(k,k)$ and $(k-1,k-1)$ it then follows from [CKS] that
$\Phi$ acts trivially $I^{k,k}$ and $I^{k-1,k-1}$.  Consequently,
$[\b^{0,0},\b^{-1,-1}] = 0$ since $\b^{-1,-1}$ maps $I^{k,k}$ to 
$I^{k-1,k-1}$ and annihilates the remaining summands appearing in 
$(9.4)$.  Thus, equation $(9.3)$ simplifies to 
$$
	\frac{d}{dy}\b^{-1,-1} =   2i[\b^{0,-1},\b^{-1,0}]	
$$
wherefrom
$$
	\b^{-1,-1} = \mu + 2i\int [\b^{0,-1},\b^{-1,0}]\,dy	\tag{9.5}
$$	

\remark{Remark} The assertion that $\Phi$ must act trivially on $Gr^W_{2k}$
and $Gr^W_{2k-2}$ is a simple consequence of the fact that $\Phi_0$ must
be a morphism of Hodge structure, and hence $\Phi_0(\xx^-)$, $\Phi_0(\xx^+)$
must be of type $(-1,1)$ and $(1,-1)$ respectively.  Therefore, the purity
of $Gr^W_{2k}$ and $Gr^W_{2k-2}$ implies that $\Phi_0$ must act trivially.
As such, the equation $-8\Phi' = Q(\Phi,\Phi)$ then implies that all of
the higher coefficients of $\Phi$ must also act trivially 
$Gr^W_{2k}$ and $Gr^W_{2k-2}$.  In particular, $N_0$ and $H$ commute with 
every element of $\lam_{(F_o,W)} = Lie_{-2}(W)$.
\endremark
\vskip 3pt

\par To continue, we now observe that by $(6.20)$ we know that if $\theta(z)$ 
was a split orbit then the associated function $h(y)$ defined by Theorem 
$(6.11)$ would be given by the formula
$$
	h(y) = e^{iyN}e^{-iyN_0}y^{-H/2} = e^{iyN_{-2}}y^{-H/2}
$$
Accordingly, when $\theta(z)$ is not split we shall write
$$
	h(y) = g(y)e^{iyN_{-2}}y^{-H/2}				\tag{9.6}
$$ 
Therefore, by equation $(8.26)$, 
$$
       g^{-1}(y)\frac{dg}{dy} 
       = y^{-H/2}.\left(-\half\Phi(\h)+\frac{H}{2y}-\Psi(e)\right)
         + i\b^{-1,-1} - iN_{-2}				\tag{9.7}
$$
Setting $\mu = N_{-2}$ it then follows from equations $(9.5)$ and $(9.7)$ that
$$
	 g^{-1}(y)\frac{dg}{dy} 
         = y^{-H/2}.\left(-\half\Phi(\h)+\frac{H}{2y}-\Psi(e)\right)
         - 2\int [\b^{0,-1},\b^{-1,0}]\,dy			\tag{9.8}
$$	
where
$$
        -2\int\, [\b^{0,-1},\b^{-1,0}]\,dy
	= y^{-2}\{\cdots\} + y^{-5/2}\{\cdots\} + \cdots	\tag{9.9}
$$ 
since $\b^{0,-1}$ and $\b^{-1,0}$ have leading order term $y^{-3/2}$.
Combining equations $(9.8)$ and $(9.9)$ with Theorem $(8.27)$ it then
follows that
$$
	g^{-1}(y)\frac{dg}{dy} = \sum_{m\geq 2}\, B_m y^{-m}
$$
Thus, just as in Corollary $(8.30)$,
$$
\aligned
	g(y) &= g(\infty)(1 + g_1 y^{-1} + g_2 y^{-2} + \cdots)		\\
	g^{-1}(y) &= (1 + f_1 y^{-1} + f_2 y^{-2} + \cdots)g^{-1}(\infty)
\endaligned
$$
where $g(\infty)$ is an arbitrary element of $\tlG$ and 
$g_n$ and $f_n$ can be expressed as universal non-commutative polynomials
in the coefficients $B_k$.

\par Continuing the analogy with \S 8, it remains therefore to show that
we can select data $(g(\infty),\{T_n\},\{S_n\})$ such that 
$$
	h(y).F_o = e^{iyN}.F
$$
In particular, the proofs of Theorem $(8.33)$ and $(8.36)$ imply 
mutatis mutandis that  $h(y).F_o = e^{iyN}.\tilde F$ where
$$
	\tilde F 
	= g(\infty)\left(
	      1+\sum_{k>0}\, \frac{1}{k!}(-i)^k(\ad\, N_0)^k\,g_k\right).
	      \hat F						\tag{9.10}
$$
provided $g(\infty)\in\ker(\ad\, N) = \ker(\ad\,N_0)\cap\ker(\ad\,N_{-2})$.
Furthermore, just as in \S 8, for purely formal algebraic reasons 
(cf\. [CKS])
$$ 
	1+\sum_{k>0}\, \frac{1}{k!}(-i)^k(\ad\, N_0)^k\,g_k
	= \exp\left(\sum_{k>0}\, Q_k(C_2,\dots,C_{k+1})\right)	 \tag{9.11}
$$
where $C_{\ell+1} = \frac{(-i)^{\ell}}{\ell!}(\ad\,N_0)^{\ell} B_{\ell+1}$.
Recycling the argument of Calculation $(8.52)$, one then finds that
$$
	C_{\ell+1}
	= i\sum_{p,q\geq 1, p+q\geq \ell+1}\, 
	  b^{\ell-1}_{p-1,q-1} \eta^{-p,-q}			\tag{9.12}
$$
where $\eta = \sum_{n>0}\, [\b_n]^{\ad\,H}_{-n}$ and 
$\b(y) = N_{-2} + \sum_{n\geq 0}\,\b_n y^{-1-n/2}$ is the series expansion of 
$\b$.

\par To complete the proof of Theorem $(4.2)$ for orbits of type $(\text{II})$,
observe that since $N_0$ acts trivially on $Gr^W_{2k}$ and $Gr^W_{2k-2}$, the
corresponding limiting mixed Hodge structure on these graded pieces is also
of type $(k,k)$ and $(k-1,k-1)$.  Therefore, if we decompose the splitting
$$
	(F,\rel W) = (e^{i\delta}.\hat F,\rel W)
$$
of the limiting mixed Hodge structure of $\theta(z)$ as
$$
	\delta = \delta_0 + \delta_{-1} + \delta_{-2}
$$
relative to the grading $Y$ defined by application of Theorem $(3.16)$
to $\hat\theta(z)$ then $\delta_0$ acts trivially on $I^{k,k}$ and 
$I^{k-1,k-1}$.  Consequently, $\delta_0$ commutes with every element
of $Lie_{-2}(W)$, and hence
$$
	e^{i\delta} = e^{i\delta_{-2}}e^{i\delta_0 + i\delta_{-1}}
$$
Proceeding as in the last part of \S 8, we can therefore pick elements
$\eta$ and $\zeta'$ so that
$$
	e^{i\delta_0 + i\delta_{-1}}	
	= e^{\zeta'}\exp\left(\sum_{k>0}\, Q_k(C_2,\dots,C_{k+1}) \right)
	\mod \exp(Lie_{-2}(W))
$$
Accordingly, since $\tlG$ contains $\exp(Lie_{-2}(W))$, we can therefore
pick elements $\eta$ and $\zeta$ such that
$$
	e^{i\delta} 
	= e^{\zeta}\left(
	      1+\sum_{k>0}\, \frac{1}{k!}(-i)^k(\ad\, N_0)^k\,g_k\right).
$$
The remaining details of Theorem $(4.2)$ regarding the uniqueness of $\eta$
and $\zeta$ etc\. are left to the reader.




\enddocument